\documentclass[secthm,seceqn,amsthm,ussrhead,12pt]{amsart}
\usepackage[utf8]{inputenc}
\usepackage[english]{babel}

\usepackage{amssymb,amsmath,amsthm,amsfonts,xcolor,enumerate,hyperref,comment,longtable,cleveref}
\usepackage{soul}
\usepackage{times}
\usepackage{cite}
\usepackage{pdflscape}
\usepackage{ulem}
\usepackage[mathcal]{euscript}
\usepackage{tikz}
\usepackage{hyperref}
\usepackage{cancel}
\usepackage{stmaryrd}

\usetikzlibrary{arrows}

\newcommand{\n}{{\bf N}}

\newtheorem{theorem}{Theorem}

\newtheorem{corollary}[theorem]{Corollary}

\theoremstyle{definition}
\newtheorem{definition}[theorem]{Definition}

\newenvironment{Proof}[1][Proof.]{\begin{trivlist}
\item[\hskip \labelsep {\bfseries #1}]}{\flushright
$\Box$\end{trivlist}}

\usepackage{stmaryrd}
\usepackage{xcolor}

\newcommand{\aut}[1]{\operatorname{\mathrm{Aut}}{#1}}

\newcommand{\cd}[2]{\mathfrak{CD}^{#1}_{#2}}

\newcommand{\la}{\langle}
\newcommand{\ra}{\rangle}
\newcommand{\La}{\Big\langle}
\newcommand{\Ra}{\Big\rangle}

\newcommand{\Dt}[2]{\Delta_{#1#2}}
\newcommand{\Dl}[2]{[\Delta_{#1#2}]}
\newcommand{\nb}[1]{\nabla_{#1}}
\newcommand{\D}[2]{\mathbf{D}^{#1}_{#2}}
\newcommand{\lb}{\lambda}
\newcommand{\0}{\theta}
\newcommand{\af}{\alpha}
\newcommand{\bt}{\beta}
\newcommand{\gm}{\gamma}
\newcommand{\dt}{\delta}
\newcommand{\eps}{\varepsilon}

\begin{document}

{\Large\noindent 
The algebraic classification of nilpotent algebras
\footnote{The work was partially  supported by  
CNPq 404649/2018-1, 302980/2019-9; RFBR 20-01-00030; by the Fundação para a Ciência e a Tecnologia (Portuguese Foundation for Science and Technology) through the project PTDC/MAT-PUR/31174/2017; by CMUP, which is financed by national funds through FCT---Funda\c c\~ao para a Ci\^encia e a Tecnologia, I.P., under the project with reference UIDB/00144/2020.
}
}

\ 

   {\bf Ivan Kaygorodov$^{a},$ Mykola Khrypchenko$^{b}$ \&   Samuel A.\ Lopes$^{c}$  \\

    \medskip
}

{\tiny

$^{a}$ CMCC, Universidade Federal do ABC, Santo Andr\'e, Brazil

$^{b}$ Departamento de Matem\'atica, Universidade Federal de Santa Catarina, Florian\'opolis, Brazil
 
$^{c}$ CMUP, Faculdade de Ci\^encias, Universidade do Porto, Rua do Campo Alegre 687, 4169-007 Porto, Portugal

\

\smallskip

   E-mail addresses:

\smallskip
    Ivan Kaygorodov (kaygorodov.ivan@gmail.com)

\smallskip  
    Mykola Khrypchenko (nskhripchenko@gmail.com)    

\smallskip    
    Samuel A.\ Lopes (slopes@fc.up.pt)

}

\

\ 

\

\noindent{\bf Abstract}: 
{\it We give the complete algebraic classification of all complex $4$-dimensional nilpotent algebras. The final list has $234$ (parametric families of) isomorphism classes of algebras, $66$ of which are new in the literature.}

\


\ 


\section*{Introduction}

The description, up to isomorphism, of all the algebras of some fixed dimension satisfying certain properties (the so-called algebraic classification) is a classical problem in algebra.
There are many papers devoted to algebraic classification of small-dimensional algebras in several varieties of
associative and non-associative algebras \cite{ack, lisa,  cfk18,         degr3, usefi1, degr1, degr2, fkkv,  kk20,   kkl19, demir, kpv19,   hac16,  ha16, maz79,   kkp19geo }.
Another interesting direction in the classification of algebras is the geometric classification (see, \cite{  fkkv,kpv19,kkp19geo     } and references therein).
Restricting our consideration to subvarieties of complex nilpotent $4$-dimensional algebras, we mention here the results on 
associative \cite{degr1},
commutative \cite{fkkv},
bicommutative \cite{kpv19},
Leibniz \cite{demir},
terminal \cite{kkp19geo} and
$\mathfrak{CD}$-algebras \cite{kk20}.
In the present paper, we complete the algebraic classification of allcomplex
$4$-dimensional nilpotent algebras. Namely, we find $66$ new algebras and parametric families of algebras completing the list of $4$-dimensional nilpotent algebras initiated in the above mentioned works.

Our approach is based on the calculation of central extensions of nilpotent algebras of dimension less than $4$. This method was developed by Skjelbred and Sund for Lie algebras in   \cite{ss78} and has been an important tool for years (see, for example, \cite{hac16, omirov, ss78} and references therein). 

In Section~\ref{SS:method} we give a brief description of the method of central extensions and its adaptation to our case. It turns out that most of the nilpotent $4$-dimensional algebras are $\frak{CD}$-algebras, which were classified in \cite{kk20}. Our classification will thus be carried out modulo this subvariety of $4$-dimensional nilpotent algebras, as explained in Section~\ref{sec-non-CD}. In Section~\ref{sec-spec-types} we discuss several other useful classes of nilpotent algebras of dimension at most $4$. In particular, we give a list of the $3$-dimensional  nilpotent  algebras (they are all, in fact, $\frak{CD}$-algebras). The remainder of Section~\ref{S:alg} (Sections~\ref{sec-CD^3_01}--\ref{sec-CD^3_04}) is devoted to the construction and classification of $1$-dimensional non-$\frak{CD}$ central extensions of $3$-dimensional nilpotent  algebras. The work done in Section~\ref{S:alg} culminates in the full classification of $4$-dimensional nilpotent algebras, which is completed in Section~\ref{S:class}, where we present the full list of representatives of isomorphism classes. Finally, in Section~\ref{S:apps} we present a few applications of our results to the classes of Lie-admissible and Alia type algebras.

\medskip
\paragraph{\bf Main result}
The main result of the paper is Theorem~\ref{teo-alg} giving a full classification of complex $4$-dimensional nilpotent algebras.
The complete list of non-isomorphic algebras consists of four parts:
\begin{enumerate}
    \item trivial $\mathfrak{CD}$-algebras, which were classified in 
    \cite[Theorem 2.1, Theorem 2.3  and Theorem 2.5]{demir}, the only anticommutative 
    trivial $\mathfrak{CD}$-algebra being $\mathfrak{CD}_{03}^{4*}$;
    
    \item terminal non-trivial extensions of the family $\cd {3*}{04}(\lambda)$, which were classified in \cite[1.4.5. 1-dimensional central extensions of ${\bf T}^{3*}_{04}$]{kkp19geo};

    \item non-trivial non-terminal $\mathfrak{ CD}$-algebras, which were classified in \cite{kk20};

    \item the nilpotent algebras found in Sections~\ref{sec-CD^3_01}--\ref{sec-CD^3_04} of the present paper.

\end{enumerate}

\section{The algebraic classification of nilpotent algebras}\label{S:alg}

\subsection{Method of classification of nilpotent algebras}\label{SS:method}
Let ${\bf A}$ be an algebra, ${\bf V}$ a vector space and ${\rm Z}^{2}\left( {\bf A},{\bf V}\right)\cong {\rm Hom}({\bf A}\otimes {\bf A},\bf V)$ denote the space of bilinear maps $\theta :{\bf A}\times 
{\bf A}\longrightarrow {\bf V}.$ For $f\in{\rm Hom}({\bf A},{\bf V})$, we introduce $\delta f\in {\rm Z}^{2}\left( {\bf A},{\bf V}\right)$ by the equality $\delta f\left( x,y\right) =f(xy)$ and
define ${\rm B}^{2}\left( {\bf A},{\bf V}\right) =\left\{\delta f \mid f\in {\rm Hom}\left( {\bf A},{\bf V}\right) \right\} $. One
can easily check that ${\rm B}^{2}({\bf A},{\bf V})$ is a linear subspace of ${\rm Z}^{2}\left( {\bf A},{\bf V}\right)$.
Let us define $\rm {H}^{2}\left( {\bf A},{\bf V}\right) $ as the quotient space ${\rm Z}^{2}\left( {\bf A},{\bf V}\right) \big/{\rm B}^{2}\left( {\bf A},{\bf V}\right)$.
The equivalence class of $\theta\in {\rm Z}^{2}\left( {\bf A%
},{\bf V}\right)$ in $\rm {H}^{2}\left( {\bf A},{\bf V}\right)$ is denoted by $\left[ \theta \right]$. 

Suppose now that $\dim{\bf A}=m<n$ and $\dim{\bf V}=n-m$. For any
bilinear map $\theta :{\bf A}\times {\bf A}\longrightarrow {\bf V%
}$, one can define on the space ${\bf A}_{\theta }:={\bf A}\oplus 
{\bf V}$ the bilinear product  $\left[ -,-\right] _{%
{\bf A}_{\theta }}$ by the equality $\left[ x+x^{\prime },y+y^{\prime }\right] _{%
{\bf A}_{\theta }}= xy  +\theta \left( x,y\right) $ for  
$x,y\in {\bf A},x^{\prime },y^{\prime }\in {\bf V}$. The algebra ${\bf A}_{\theta }$ is called an $(n-m)$-{\it dimensional central extension} of ${\bf A}$ by ${\bf V}$.
It is also clear that ${\bf A}_{\theta }$ is nilpotent if and only if  ${\bf A}$ is.

For a bilinear map $\theta :{\bf A}\times {\bf A}\longrightarrow {\bf V}$, the space $\theta ^{\bot }=\left\{ x\in {\bf A}\mid \theta \left(
x,{\bf A}\right) =\theta \left(
{\bf A},x\right) =0\right\} $ is called the {\it annihilator} of $\theta$.
For an algebra ${\bf A}$, the ideal 
${\rm Ann}\left( {\bf A}\right) =\left\{ x\in {\bf A}\mid x{\bf A} ={\bf A}x =0\right\}$ is called the {\it annihilator} of ${\bf A}$.
One has
\begin{equation*}
{\rm Ann}\left( {\bf A}_{\theta }\right) =\left( \theta ^{\bot }\cap {\rm Ann}\left( 
{\bf A}\right) \right) \oplus {\bf V}.
\end{equation*}
Any $n$-dimensional algebra with non-trivial annihilator can be represented in
the form ${\bf A}_{\theta }$ for some $m$-dimensional algebra ${\bf A}$, an $(n-m)$-dimensional vector space ${\bf V}$ and $\theta \in {\rm Z}^{2}\left( {\bf A},{\bf V}\right)$, where $m<n$ (see \cite[Lemma 5]{hac16}).
Moreover, there is a unique such representation with $m=n-\dim{\rm Ann}({\bf A})$. Note also that the latter equality is equivalent to the condition  $\theta ^{\bot }\cap {\rm Ann}\left( 
{\bf A}\right)=0$. 

Let us pick some $\phi\in {\rm Aut}\left( {\bf A}\right)$, where ${\rm Aut}\left( {\bf A}\right)$ is the automorphism group of  ${\bf A}$.
For $\theta\in {\rm Z}^{2}\left( {\bf A},{\bf V}\right)$, let us define $(\phi \theta) \left( x,y\right) =\theta \left( \phi \left( x\right)
,\phi \left( y\right) \right) $. Then we get an action of ${\rm Aut}\left( {\bf A}\right) $ on ${\rm Z}^{2}\left( {\bf A},{\bf V}\right)$ that induces an action of the same group on $\rm {H}^{2}\left( {\bf A},{\bf V}\right)$.  

\begin{definition}
Let ${\bf A}$ be an algebra and $I$ be a subspace of ${\rm Ann}({\bf A})$. If ${\bf A}={\bf A}_0 \oplus I$, for some subalgebra ${\bf A}_0$ of ${\bf A}$,
then $I$ is called an {\it annihilator component} of ${\bf A}$. We say that ${\bf A}$ is {\it split} if it has some non-zero annihilator component; otherwise we say that ${\bf A}$ is {\it non-split}.
\end{definition}

For a linear space $\bf U$, the {\it Grassmannian} $G_{s}\left( {\bf U}\right) $ is
the set of all $s$-dimensional linear subspaces of ${\bf U}$. For any $s\ge 1$, the action of ${\rm Aut}\left( {\bf A}\right)$ on $\rm {H}^{2}\left( {\bf A},\mathbb{C}\right)$ induces 
an action of the same group on $G_{s}\left( \rm {H}^{2}\left( {\bf A},\mathbb{C}\right) \right)$.
Let us define
$$
{\bf T}_{s}\left( {\bf A}\right) =\left\{ {\bf W}\in G_{s}\left( \rm {H}^{2}\left( {\bf A},\mathbb{C}\right) \right)\left|\underset{[\theta]\in W}{\cap }\theta^{\bot }\cap {\rm Ann}\left( {\bf A}\right) =0\right.\right\}.
$$
Note that ${\bf T}_{s}\left( {\bf A}\right)$ is stable under the action of ${\rm Aut}\left( {\bf A}\right) $.

Let us fix a basis $e_{1},\ldots,e_{s} $ of ${\bf V}$ and $\theta \in {\rm Z}^{2}\left( {\bf A},{\bf V}\right) $. Then there are unique $\theta _{i}\in {\rm Z}^{2}\left( {\bf A},\mathbb{C}\right)$ ($1\le i\le s$) such that $\theta \left( x,y\right) =\underset{i=1}{\overset{s}{%
\sum }}\theta _{i}\left( x,y\right) e_{i}$ for all $x,y\in{\bf A}$. 
Note that $\theta ^{\bot
}=\theta^{\bot} _{1}\cap \theta^{\bot} _{2}\cdots \cap \theta^{\bot} _{s}$ in this case.
If   $\theta ^{\bot
}\cap {\rm Ann}\left( {\bf A}\right) =0$, then by \cite[Lemma 13]{hac16} the algebra ${\bf A}_{\theta }$ is split if and only if $\left[ \theta _{1}\right],\left[\theta _{2}\right] ,\ldots ,\left[ \theta _{s}\right] $ are linearly
dependent in $\rm {H}^{2}\left( {\bf A},\mathbb{C}\right)$. Thus, if $\theta ^{\bot
}\cap {\rm Ann}\left( {\bf A}\right) =0$ and ${\bf A}_{\theta }$ is non-split, then $\left\langle \left[ \theta _{1}\right] , \ldots,%
\left[ \theta _{s}\right] \right\rangle$ is an element of ${\bf T}_{s}\left( {\bf A}\right)$.
Now, if $\vartheta\in {\rm Z}^{2}\left( {\bf A},\bf{V}\right)$ is such that $\vartheta ^{\bot
}\cap {\rm Ann}\left( {\bf A}\right) =0$ and ${\bf A}_{\vartheta }$ is non-split, then by \cite[Lemma 17]{hac16} one has ${\bf A}_{\vartheta }\cong{\bf A}_{\theta }$ if and only if
$\left\langle \left[ \theta _{1}\right] ,\left[ \theta _{2}%
\right] ,\ldots ,\left[ \theta _{s}\right] \right\rangle,
\left\langle \left[ \vartheta _{1}\right] ,\left[ \vartheta _{2}\right] ,\ldots,%
\left[ \vartheta _{s}\right] \right\rangle\in {\bf T}_{s}\left( {\bf A}\right)$ belong to the same orbit under the action of ${\rm Aut}\left( {\bf A}\right) $, where $%
\vartheta \left( x,y\right) =\underset{i=1}{\overset{s}{\sum }}\vartheta
_{i}\left( x,y\right) e_{i}$.

Hence, there is a one-to-one correspondence between the set of $%
{\rm Aut}\left( {\bf A}\right) $-orbits on ${\bf T}_{s}\left( {\bf A}%
\right) $ and the set of isomorphism classes of central extensions of $\bf{A}$ by $\bf{V}$ with $s$-dimensional annihilator and trivial annihilator component.
Consequently, to construct all $s$-dimensional central extensions with trivial annihilator component
of a given $(n-s)$-dimensional algebra ${\bf A}$ one has to describe ${\bf T}_{s}({\bf A})$, ${\rm Aut}({\bf A})$ and the action of ${\rm Aut}({\bf A})$ on ${\bf T}_{s}({\bf A})$ and then
for each orbit under the action of ${\rm Aut}({\bf A})$ on ${\bf T}_{s}({\bf A})$ pick a representative and construct the algebra corresponding to it.

\subsubsection{Reduction to non-$\mathfrak{CD}$-algebras}\label{sec-non-CD}
The class of  $\mathfrak{CD}$-algebras is defined by the 
property that the commutator of any pair of multiplication operators is a derivation \cite{ack, kk20}; 
namely, an algebra $\mathfrak{A}$ is a $\mathfrak{CD}$-algebra if and only if  
\[ [T_x,T_y]   \in \mathfrak{Der} (\mathfrak{A}),\]
for all $x,y \in \mathfrak{A}$, where  $T_z \in \{ R_z,L_z\}$. Here we use the notation $R_z$ (resp. $L_z$) for the operator of right (resp. left) multiplication in $\mathfrak{A}$. We will denote the variety of $\mathfrak{CD}$-algebras by $\mathfrak{CD}$.
In terms of identities, the class of $\mathfrak{CD}$-algebras is defined by the following three:
\begin{align*}
    ((xy)a)b-((xy)b)a&=((xa)b-(xb)a)y+x((ya)b-(yb)a),\\
    (a(xy))b-a((xy)b)&=((ax)b-a(xb))y+x((ay)b-a(yb)),\\
    a(b(xy))-b(a(xy))&=(a(bx)-b(ax))y+x(a(by)-b(ay)).
\end{align*}

Our method of classification of nilpotent algebras will be based on a classification of 
 $\mathfrak{CD}$-algebras~\cite{kk20}. Clearly, any central extension of a non-$\mathfrak{CD}$-algebra is a non-$\mathfrak{CD}$-algebra. But a $\mathfrak{CD}$-algebra may have extensions which are not $\mathfrak{CD}$-algebras. More precisely, let $\bf{A}$ be a $\mathfrak{CD}$-algebra and $\theta \in {\rm Z^2}(\bf{A}, {\mathbb C}).$ Then ${\bf{A}}_{\theta }$ is a $\mathfrak{CD}$-algebra if and only if
\begin{align}
\label{cd1}    \theta((xy)a,b)-\theta((xy)b,a)&=\theta((xa)b-(xb)a,y)+\theta(x,(ya)b-(yb)a),\\
\label{cd2}    \theta(a(xy),b)-\theta(a,(xy)b)&=\theta((ax)b-a(xb),y)+\theta(x,(ay)b-a(yb)),\\
\label{cd3}    \theta(a,b(xy))-\theta(b,a(xy))&=\theta(a(bx)-b(ax),y)+\theta(x,a(by)-b(ay)).
\end{align}
for all $x,y,z\in \bf{A}.$ 
Define the subspace ${\rm Z_\mathfrak{CD}^2}(\bf{A},{\mathbb C})$ of ${\rm Z^2}(\bf{A},{\mathbb C})$ by
\begin{equation*}
{\rm Z_{\mathfrak{CD}}^2}(\bf{A},{\mathbb C}) =\left\{\begin{array}{ll} \theta \in {\rm Z^2}(\bf{A},{\mathbb C}): &  
\theta \mbox{ satisfies }
(\ref{cd1}), (\ref{cd2}) \mbox{ and } (\ref{cd3}) \end{array}\right\}.
\end{equation*}

Observe that ${\rm B^2}(\bf{A},{\mathbb C})\subseteq{\rm Z_\mathfrak{CD}^2}(\bf{A},{\mathbb C}).$
Let ${\rm H_\mathfrak{CD}^2}(\bf{A},{\mathbb C}) =%
{\rm Z_\mathfrak{CD}^2}(\bf{A},{\mathbb C}) \big/{\rm B^2}(\bf{A},{\mathbb C}).$ Then ${\rm H_\mathfrak{CD}^2}(\bf{A},{\mathbb C})$ is a subspace of $%
{\rm H^2}(\bf{A},{\mathbb C}).$ Define
\[{\bf R}_{s}(\bf{A}) =\left\{ {\bf W}\in {\bf T}_{s}(\bf{A}) :{\bf W}\in G_{s}({\rm H_\mathfrak{CD}^2}(\bf{A},{\mathbb C}) ) \right\}, \]
\[	{\bf U}_{s}(\bf{A}) =\left\{ {\bf W}\in {\bf T}_{s}(\bf{A}) :{\bf W}\notin G_{s}({\rm H_\mathfrak{CD}^2}(\bf{A},{\mathbb C}) ) \right\}.\]
Then ${\bf T}_{s}(\bf{A}) ={\bf R}_{s}(\bf{A})$ $\mathbin{\mathaccent\cdot\cup}$ ${\bf U}_{s}(\bf{A}).$ The sets ${\bf R}_{s}(\bf{A}) $ and ${\bf U}_{s}(\bf{A})$ are stable under the action
of $\operatorname{Aut}(\bf{A}).$ 
Thus, the  nilpotent algebras
corresponding to the representatives of $\operatorname{Aut}(\bf{A})$-orbits on ${\bf R}_{s}(\bf{A})$ are $\mathfrak{CD}$-algebras,
while those corresponding to the representatives of $\operatorname{Aut}(\bf{A}) $-orbits on ${\bf U}_{s}(\bf{A})$ are non-$\mathfrak{CD}$-algebras. Hence, we may construct all nilpotent non-split non-$\mathfrak{CD}$-algebras $\bf{A}$ of dimension $n$ with $s$-dimensional annihilator
from a given nilpotent algebra $\bf{A}^{\prime }$ of dimension $n-s$ as follows:

\begin{enumerate}
\item If $\bf{A}^{\prime }$ is non-$\mathfrak{CD}$, then apply the procedure.

\item Otherwise, do the following:

\begin{enumerate}
\item Determine ${\bf U}_{s}(\bf{A}^{\prime})$ and $\operatorname{Aut}(\bf{A}^{\prime }).$

\item Determine the set of $\operatorname{Aut}(\bf{A}^{\prime })$-orbits on ${\bf U}_{s}(\bf{A}^{\prime }).$

\item For each orbit, construct the nilpotent algebra corresponding to one of its
representatives.
\end{enumerate}
\end{enumerate}

We will use the following auxiliary notation during the construction of central extensions.
Let ${\bf A}$ be an algebra with basis $e_{1},e_{2},\ldots,e_{n}$.  Then $\Delta_{ij}:{\bf A}\times {\bf A}\longrightarrow \mathbb{C}$ denotes the   bilinear form  defined by the equalities $\Delta _{ij}\left( e_{i},e_{j}\right) =1$
and $\Delta _{ij}\left( e_{l},e_{m}\right) =0$ for $(l,m) \neq (i,j)$. In this case, the $\Delta_{ij}$ with $1\leq i , j\leq n $ form a basis of the space of  bilinear forms on $\bf{A}$. 

We also denote by

\begin{longtable}{lll}
$\mathfrak{CD}^{i*}_j$& the $j$th $i$-dimensional nilpotent trivial  $\mathfrak{CD}$-algebra, \\
$\mathfrak{CD}^i_j$&  the $j$th $i$-dimensional nilpotent  non-trivial  $\mathfrak{CD}$-algebra, \\
$\D{i}{j}$&  the $j$th $i$-dimensional nilpotent  terminal  algebra,\\
${\n}^i_j$& the $j$th $i$-dimensional nilpotent non-$\mathfrak{CD}$-algebra.
\end{longtable}


\subsection{Some special types of nilpotent algebras}\label{sec-spec-types}
\subsubsection{Trivial $\mathfrak{CD}$-algebras}
Recall that the class of $n$-dimensional algebras defined by the identities $(xy)z=0$ and $x(yz)=0$ 
lies in the intersection of all well-known varieties of algebras defined by a family of polynomial identities of degree $3,$ such as Leibniz, Zinbiel, associative, Novikov and many other algebras.
On the other hand, 
all algebras defined by the identities $(xy)z=0$ and $x(yz)=0$ are central extensions of some suitable algebra with zero product.
The list of all  non-anticommutative $4$-dimensional algebras defined by the identities $(xy)z=0$ and $x(yz)=0$  can be found in \cite{demir}.
Note that there is only one $4$-dimensional nilpotent anticommutative algebra satisfying the identity  $(xy)z=0.$
Obviously, all algebras from this list are  $4$-dimensional nilpotent $\mathfrak{CD}$-algebras.

\subsubsection{$2$-dimensional nilpotent algebras}
There is only one non-zero $2$-dimensional nilpotent algebra.
It is a $\mathfrak{CD}$-algebra:

\begin{longtable}{ll llll}
$\cd {2*}{01}$ & $:$ & $e_1 e_1 = e_2$
\end{longtable}

\subsubsection{$3$-dimensional nilpotent algebras}
Thanks to \cite{cfk18}, we have the classification of all $3$-dimensional nilpotent algebras.
It is easy to see that each $3$-dimensional nilpotent algebra is a $\mathfrak{CD}$-algebra:

\begin{longtable}{lllll llll}
$\cd 3{01}$&$:$& $e_1 e_1 = e_2$ & $e_2 e_2=e_3$ \\
\hline
$\cd 3{02}$&$:$& $e_1 e_1 = e_2$ & $e_2 e_1= e_3$ & $e_2 e_2=e_3$ \\
\hline$\cd 3{03}$&$:$& $e_1 e_1 = e_2$ & $e_2 e_1=e_3$ \\
\hline$\cd 3{04}(\lambda)$&$:$& $ e_1 e_1 = e_2$ & $e_1 e_2=e_3$ & $e_2 e_1=\lambda e_3$ \\
\hline$\cd {3*}{01}$&$:$& $e_1 e_1 = e_2$\\
\hline$\cd {3*}{02}$&$:$& $e_1 e_1 = e_3$ &$ e_2 e_2=e_3$ \\
\hline$\cd {3*}{03}$&$:$& $e_1 e_2=e_3$ & $e_2 e_1=-e_3$ \\
\hline$\cd {3*}{04}(\lambda)$&$:$& $e_1 e_1 = \lambda e_3$ & $e_2 e_1=e_3$  & $e_2 e_2=e_3$
\end{longtable}

\subsubsection{$4$-dimensional nilpotent algebras with $2$-dimensional annihilator }
Thanks to \cite{cfk18}, we have the classification of all $4$-dimensional non-split nilpotent non-trivial algebras with $2$-dimensional annihilator.
All of these algebras are $\mathfrak{CD}$-algebras:
\begin{longtable}{lllll llll}
$\cd 4{05}$&:  & $e_1 e_1 = e_2$ & $e_2 e_1=e_4$ & $e_2 e_2=e_3$ \\
\hline$\cd 4{06}$ & $:$ & $ e_1 e_1 = e_2$ & $e_1 e_2=e_4$ & $e_2 e_1=e_3$  \\
\hline$\cd 4{07}(\lambda)$& $:$ & $e_1 e_1 = e_2$ & $e_1 e_2=e_4$ & $e_2 e_1=\lambda e_4$ & $e_2 e_2=e_3$ \\
\end{longtable}

\subsubsection{$4$-dimensional nilpotent algebras with $1$-dimensional annihilator }
Thanks to \cite{kk20} we have the classification of all $4$-dimensional nilpotent $\mathfrak{CD}$-algebras with $1$-dimensional annihilator.
The remaining nilpotent non-$\mathfrak{CD}$-algebras with $1$-dimensional annihilator will be found in the present paper.

\subsubsection{Trivial central extensions.}
Thanks to \cite{kk20} we know that all central extensions of 
$\mathfrak{CD}^{3*}_{01},$
$\mathfrak{CD}^{3*}_{02},$
$\mathfrak{CD}^{3*}_{03}$ and
$\mathfrak{CD}^{3*}_{04}$  are $\mathfrak{CD}$-algebras.

\bigskip\bigskip
 
\subsection{$1$-dimensional central extensions of $\cd 3{01}$}\label{sec-CD^3_01}
Here we collect all information about ${\mathfrak{CD}_{01}^{3}}$:

\begin{longtable}{|l|l|l|}
\hline
Algebra  & Multiplication & Cohomology   \\
\hline
${\mathfrak{CD}}_{01}^{3}$ &  
$\begin{array}{l}e_1 e_1 = e_2\\  e_2 e_2=e_3
\end{array}$& 
$\begin{array}{lcl}
{\rm H^2_{\mathfrak{CD}}}({\mathfrak{CD}}_{01}^{3})&=&\langle  [\Delta_{12}], [\Delta_{21}] \rangle \\
{\rm H^2}({\mathfrak{CD}}_{01}^{3})&=&{\rm H^2_{\mathfrak{CD}}}({\mathfrak{CD}_{01}^3}) \oplus 
\langle[\Delta_{13}], [\Delta_{31}], [\Delta_{23}], [\Delta_{32}],  [\Delta_{33}]    \rangle 
\end{array}$\\
\hline
\end{longtable}	
	Let us use the following notations 
	\begin{longtable}{llll}
$\nb 1 = \Dl 12$&$ \nb 2 = \Dl 21$&  $\nb 3 = \Dl 13$ &$\nb 4=\Dl 31$\\
$\nb 5 = \Dl 23$&$ \nb 6 = \Dl 32$&  $\nb 7 = \Dl 33$. 

	\end{longtable}
	Take $\0=\sum_{i=1}^7\af_i\nb i\in {\rm H}^2  (\cd 3{01})$.
	If 
	$$
	\phi=
	\begin{pmatrix}
	x &    0  &  0\\
	0 &  x^2  &  0\\
	y &   0  &  x^4
	\end{pmatrix}\in\aut{\cd 3{01}},
	$$
	then
	$$
	\phi^T\begin{pmatrix}
	0        &  \af_1  & \af_3\\
	\af_2        &  0 & \af_5\\
	\af_4  &     \af_6   & \af_7
	\end{pmatrix} \phi=
	\begin{pmatrix}
	\af^*    &  \af_1^* & \af_3^*  \\
	\af^*_2  &  0      & \af^*_5\\
	\af^*_4  &  \af^*_6  & \af^*_7
	\end{pmatrix},
	$$
	where

\begin{longtable}{llll}
$\af^*_1=x^2 (x \alpha_1+y \alpha_6)$ & 
$\af^*_2=x^2 (x \alpha_2+y \alpha_5)$ & 
$\af^*_3=x^4 (x \alpha_3+y \alpha_7)$ & 
$\af^*_4=x^4 (x \alpha_4+y \alpha_7)$\\
$\af^*_5=x^6 \alpha_5$ &
$\af^*_6=x^6 \alpha_6$& 
$\af^*_7=x^8 \alpha_7$.&
	\end{longtable} 	
 Hence, $\phi\langle\0\rangle=\langle\0^*\rangle$, where $\0^*=\sum\limits_{i=1}^7 \af_i^*  \nb i.$
We are only interested in elements with $(\alpha_3, \af_4, \af_5, \af_6, \af_7)\neq (0,0,0,0,0).$
Then

\begin{enumerate}

\item $\af_7\neq 0.$ 
    \begin{enumerate} 
    \item $\af_6\neq 0$.
    Then choosing $x=\sqrt{\af_6 \af_7^{-1}}$ and $y=- \af_4\af_7^{-1}x$, we have the family of representatives 
    $\langle \af  \nabla_1+\bt \nabla_2+\gm \nabla_3+\dt \nb 5+\nb 6+\nb 7 \rangle.$
Observe that two distinct quadruples $(\af,\bt,\gm, \dt)$ and $(\af',\bt',\gm', \dt')$ determine the same orbit if and only if $(\af,\bt,\gm, \dt)=(-\af',-\bt',-\gm', \dt').$   
   
    \item  $\af_6= 0$ and $\af_5\neq0$.
    Choosing $x=\sqrt{\af_5 \af_7^{-1}}$ and $y=- \af_4\af_7^{-1}x$, we have the family of representatives 
    $\langle \af  \nabla_1+\bt \nabla_2+\gm \nabla_3+  \nb 5 +\nb 7 \rangle.$
Observe that two distinct triples $(\af,\bt,\gm)$ and $(\af',\bt',\gm')$ determine the same orbit if and only if  $(\af,\bt,\gm)=(-\af',-\bt',-\gm').$   
       
    \item  $\af_6=\af_5=0$ and $\af_3 \neq \af_4$.
    Choosing $x=\sqrt[3]{(\af_3-\af_4) \af_7^{-1}}$ and $y=- \af_4\af_7^{-1}x$, we have the family of representatives 
    $\langle \af  \nabla_1+\bt \nabla_2+  \nabla_3+   \nb 7 \rangle.$
Observe that two pairs $(\af,\bt)$ and $(\af',\bt')$ determine the same orbit if and only if $(\af,\bt )=(\xi \af',\xi \bt'),$ where $\xi^3=1.$   
    
    \item $\af_6=\af_5= 0,$ $\af_3 = \af_4$ and $\af_2\ne 0$.
    Choosing 
    $x=\sqrt[5]{\af_2\af_7^{-1}}$ and $y=-\af_4\af_7^{-1}x$, we have the family of  representatives of distinct orbits
    $\langle \af  \nabla_1+\nabla_2+\nb 7 \rangle.$
  
    \item $\af_6=\af_5=0$, $\af_3 = \af_4$ and  $\af_2=0$.
    Choosing 
     $y=-\af_4\af_7^{-1}x$, we have two representatives 
     $\langle \nb 7 \rangle$ and $\langle  \nabla_1+ \nb 7 \rangle$
    depending on whether $\af_1=0$ or not.

\end{enumerate}

\item $\af_7=0.$

\begin{enumerate}
    \item $\af_6\neq 0$ and $\af_4 \neq 0$.
    Choosing 
    $x=\af_4 \af_6^{-1}$ and $y=-\af_1\af_4 \af_6^{-2},$ we have the family of  representatives of distinct orbits
       $\langle \af  \nabla_2+\bt \nabla_3+ \nabla_4+ \gm \nb 5 +\nb 6 \rangle.$

     \item $\af_6\neq 0,$  $\af_4 = 0$ and $\af_3\neq0$.
    Choosing 
    $x=\af_3 \af_6^{-1}$ and $y=-\af_1\af_3 \af_6^{-2},$ we have the family of representatives of distinct orbits
       $\langle \af  \nabla_2+ \nabla_3+ \bt \nb 5 +\nb 6 \rangle.$

     \item $\af_6\neq 0$ and  $\af_4=\af_3 = 0$.
    Choosing 
    $y=-\af_1\af_6^{-1}x$ we have two families of representatives of distinct orbits
       $\langle  \af \nb 5 +\nb 6 \rangle$ and $\langle  \nabla_2+   \af \nb 5 +\nb 6 \rangle$
        depending on whether $\af_1\af_5 - \af_2\af_6=0$ or not.  

    \item $\af_6 = 0,$ $\af_5\neq0$ and   $\af_4 \neq 0$. 
    Choosing 
    $x=\af_4 \af_5^{-1}$ and $y=-\af_2\af_4  \af_5^{-2},$ we have the family  of representatives of distinct orbits
       $\langle  \af \nb 1+\bt \nabla_3 +\nb 4+\nb 5 \rangle.$

 \item $\af_6 = 0,$ $\af_5\neq0$ and    $\af_4 = 0$.
     Choosing $y=-\af_2\af_5^{-1}x$, we have the following representatives of distinct orbits
       $\langle  \af \nb 1 +\nb 3+\nb 5 \rangle$,
       $\langle   \nb 1 +\nb 5 \rangle$ and 
              $\langle  \nb 5 \rangle,$ corresponding to the 3 cases: $\af_3\ne 0$; $\af_3=0$, $\af_1\ne 0$; $\af_3=\af_1=0$, respectively.

 \item $\af_6 = \af_5 = 0$ and    $\af_4 \neq 0$. 
     We have the following families of representatives of distinct orbits
       $\langle   \nb 1+\af \nb 2+\bt\nb 3+\nb 4 \rangle$, $\langle \nb 2+\af\nb 3+\nb 4 \rangle$ and $\langle \af\nb 3+\nb 4 \rangle$ corresponding to the 3 cases: $\af_1\ne 0$; $\af_1=0$, $\af_2\ne 0$; $\af_1=\af_2=0$, respectively.

 \item $\af_6 = \af_5=\af_4 = 0$ and $\af_3\neq 0$.  
     We have the following representatives of distinct orbits
       $\langle  \af \nb 1+ \nabla_2 +\nb 3 \rangle,$
       $\langle   \nb 1+  \nb 3 \rangle$ and 
              $\langle  \nb 3 \rangle$ corresponding to the 3 cases: $\af_2\ne 0$; $\af_2=0$, $\af_1\ne 0$; $\af_2=\af_1=0$, respectively.

\end{enumerate}
\end{enumerate}
    Summarizing, we have the following distinct orbits:
\begin{longtable}{lll}

$\langle  \af \nb 1+ \nabla_2 +\nb 3 \rangle$ & 

\multicolumn{2}{l}{$\langle \af  \nabla_1+\bt \nabla_2+\gm \nabla_3+\dt \nb 5+\nb 6+\nb 7 \rangle^{O(\af,\bt,\gm, \dt)=O(-\af,-\bt,-\gm, \dt)}$}\\

$\langle  \af \nb 1+\nb 3+\nb 5 \rangle$ &
\multicolumn{2}{l}{$\langle \af  \nabla_1+\bt \nabla_2+\gm \nabla_3+\nb 5+\nb 7 \rangle^{O(\af,\bt,\gm)=O(-\af,-\bt,-\gm)}$}\\

\multicolumn{2}{l}{$\langle \af  \nabla_1+\bt \nabla_2+ \nabla_3 +\nb 7 \rangle^{O(\af,\bt )=O(\sqrt[3]1(\af,\bt ))}$}&


  $\langle \af  \nabla_1+\nabla_2+     \nb 7 \rangle$ \\

$\langle  \nb 1+  \nb 3 \rangle$  &

$\langle   \nb 1+\af \nb 2+\bt\nb 3+\nb 4 \rangle$  &

$\langle \af \nb 1+\bt \nabla_3 +\nb 4+\nb 5 \rangle$ \\

$\langle   \nb 1 +\nb 5 \rangle$ &

$\langle  \nabla_1+ \nb 7 \rangle$ &   

$\langle \nb 2+ \af  \nabla_3  +\nb 4 \rangle$ \\
$\langle \af  \nabla_2+\bt \nabla_3+ \nabla_4+ \gm \nb 5 +\nb 6 \rangle$ &
$\langle \af  \nabla_2+ \nabla_3+ \bt \nb 5 +\nb 6 \rangle$ &

$\langle \nb 2 + \af\nb 5+\nb 6 \rangle$  \\
$\langle \af  \nabla_3  +\nb 4 \rangle$ &
$\langle \nb 3 \rangle$ &
$\langle  \nb 5 \rangle$ \\

$\langle \af \nb 5 +\nb 6 \rangle$ &

$\langle \nb 7 \rangle$ 
\end{longtable}
They correspond to the following new algebras:

{\tiny 
\begin{longtable}{lllllllllllllllllll}
$\n^4_{01}(\af)$ &$:$& $e_1 e_1 = e_2$& $e_1e_2=\af e_4$& $e_1e_3=e_4$& $e_2e_1=e_4$& $e_2 e_2=e_3$ \\ 
\hline
$\n^4_{02}(\af,\bt,\gm,\delta)$ &$:$& 
$e_1 e_1 = e_2$& $e_1e_2=\af e_4$& $e_1e_3=\gm e_4$&  $e_2e_1=\bt e_4$\\&& 
$e_2 e_2=e_3$& $e_2e_3=\delta e_4$& $e_3e_2=e_4$& $e_3e_3=e_4$\\ 
\hline
$\n^4_{03}(\af)$ &$:$& $e_1 e_1 = e_2$& $e_1e_2=\af e_4$& $e_1e_3=e_4$& $e_2 e_2=e_3$ & $e_2e_3=e_4$\\  \hline

$\n^4_{04}(\af,\bt,\gm)$ &$:$& 
$e_1 e_1 = e_2$&$e_1e_2=\af e_4$&$e_1e_3=\gm e_4$&$e_2e_1=\bt e_4$\\
&&  $e_2 e_2=e_3$&$e_2e_3=e_4$&$e_3e_3=e_4$ \\  \hline

$\n^4_{05}(\af,\bt)$ &$:$& 
$e_1 e_1 = e_2$&$e_1e_2=\af e_4$&$e_1e_3=e_4$ &$e_2e_1=\bt e_4$&  $e_2 e_2=e_3$&$e_3e_3=e_4$ \\  \hline

$\n^4_{06}(\af)$ &$:$& $e_1 e_1 = e_2$&$e_1e_2=\af e_4$&$e_2e_1=e_4$&  $e_2 e_2=e_3$&$e_3e_3=e_4$ \\  \hline

$\n^4_{07}$ &$:$& $e_1 e_1 = e_2$&$e_1e_2=e_4$&$e_1e_3=e_4$&  $e_2 e_2=e_3$ \\  \hline

$\n^4_{08}(\af,\bt)$ &$:$& $e_1 e_1 = e_2$&$e_1e_2=e_4$&$e_1e_3=\bt e_4$&$e_2e_1=\af e_4$&  $e_2 e_2=e_3$ &$e_3e_1=e_4$\\  \hline

$\n^4_{09}(\af,\bt)$ &$:$& 
$e_1 e_1 = e_2$&$e_1e_2=\af e_4$&$e_1e_3=\bt e_4$&  $e_2 e_2=e_3$ &$e_2e_3=e_4$&$e_3e_1=e_4$ \\  \hline

$\n^4_{10}$ &$:$& $e_1 e_1 = e_2$&$e_1e_2=e_4$&  $e_2 e_2=e_3$ &$e_2e_3=e_4$ \\  \hline

$\n^4_{11}$ &$:$& $e_1 e_1 = e_2$&$e_1e_2=e_4$&  $e_2 e_2=e_3$&$e_3e_3=e_4$ \\  \hline

$\n^4_{12}(\af)$ &$:$& $e_1 e_1 = e_2$&  $e_1e_3=\af e_4$&$e_2e_1=e_4$ &  $e_2 e_2=e_3$ &$e_3e_1=e_4$ \\  \hline

$\n^4_{13}(\af,\bt,\gm)$ &$:$& 
$e_1 e_1 = e_2$&$e_1e_3=\bt e_4$&$e_2e_1=\af e_4$&  $e_2 e_2=e_3$\\
&&$e_2e_3=\gm e_4$&$e_3e_1=e_4$&$e_3e_2=e_4$ \\  \hline

$\n^4_{14}(\af, \bt)$ &$:$& $e_1 e_1 = e_2$&$e_1e_3 =e_4$&$e_2e_1= \af e_4$&  $e_2 e_2=e_3$ &$e_2e_3=\bt e_4$&$e_3e_2=e_4$ \\  \hline

$\n^4_{15}(\af)$ &$:$& $e_1 e_1 = e_2$&$e_2e_1=e_4$&  $e_2 e_2=e_3$ &$e_2e_3=\af e_4$ &$e_3e_2=e_4$ \\  \hline

$\n^4_{16}(\af)$ &$:$& $e_1 e_1 = e_2$&$e_1e_3=\af e_4$ &  $e_2 e_2=e_3$ &$e_3e_1=e_4$\\  \hline

$\n^4_{17}$ &$:$& $e_1 e_1 = e_2$&  $e_1e_3=e_4$ & $e_2 e_2=e_3$\\  \hline

$\n^4_{18}$ &$:$& $e_1 e_1 = e_2$&  $e_2 e_2=e_3$ &$e_2e_3=e_4$\\  \hline

$\n^4_{19}(\af)$ &$:$& $e_1 e_1 = e_2$&  $e_2 e_2=e_3$ &$e_2e_3=\af e_4$ &$e_3e_2=e_4$ \\  \hline

$\n^4_{20}$ &$:$& $e_1 e_1 = e_2$&  $e_2 e_2=e_3$ &$e_3e_3=e_4$\\  \hline

\end{longtable}
}
 \begin{longtable}{cc}
$\n^4_{02}(\af,\bt,\gm,\delta) \cong \n^4_{02}(-(\af,\bt,\gm),\delta)$ &
$\n^4_{04}(\af,\bt,\gm) \cong \n^4_{04}(-(\af,\bt,\gm))$ \\
\multicolumn{2}{c}{$\n^4_{05}(\af,\bt) \cong \n^4_{05}(\sqrt[3]{1}(\af,\bt))$}  
 \end{longtable}

 	\subsection{$1$-dimensional central extensions of $\cd 3{02}$}
Here we collect all information about ${\mathfrak{CD}_{02}^{3}}$:
\begin{longtable}{|l|l|l|}
\hline
Algebra  & Multiplication & Cohomology   \\
\hline
${\mathfrak{CD}}_{02}^{3}$ &  
$\begin{array}{l}
e_1 e_1 = e_2\\ 
e_2 e_1= e_3\\ 
e_2 e_2=e_3
\end{array}$& 
$\begin{array}{lcl}
{\rm H^2_{\mathfrak{CD}}}({\mathfrak{CD}}_{02}^{3})&=&\langle  [\Delta_{12}], [\Delta_{21}] \rangle \\
{\rm H^2}({\mathfrak{CD}}_{02}^{3})&=&{\rm H^2_{\mathfrak{CD}}}({\mathfrak{CD}_{02}^3}) \oplus 
\langle[\Delta_{13}], [\Delta_{31}], [\Delta_{23}], [\Delta_{32}],  [\Delta_{33}]    \rangle 
\end{array}$\\
\hline
\end{longtable}		
	Let us use the following notations 
	\begin{longtable}{llll}
$\nb 1 = \Dl 12$&$ \nb 2 = \Dl 21$&  $\nb 3 = \Dl 13$ &$\nb 4=\Dl 31$\\
$\nb 5 = \Dl 23$&$ \nb 6 = \Dl 32$&  $\nb 7 = \Dl 33$. 

	\end{longtable}
	Take $\0=\sum_{i=1}^7\af_i\nb i\in {\rm H}^2(\cd 3{02})$.
	If 
	$$
	\phi=
	\begin{pmatrix}
	1 &    0  &  0\\
	0 & 1  &  0\\
	x &   0  & 1
	\end{pmatrix}\in\aut{\cd 3{02}},
	$$
	then
	$$
	\phi^T\begin{pmatrix}
	0        &  \af_1  & \af_3\\
	\af_2        &  0 & \af_5\\
	\af_4  &     \af_6   & \af_7
	\end{pmatrix} \phi=
	\begin{pmatrix}
	\af^*    &  \af_1^* & \af_3^*  \\
	\af^*_2  &  0      & \af^*_5\\
	\af^*_4  &  \af^*_6  & \af^*_7
	\end{pmatrix},
	$$
	where

\begin{longtable}{llll}
$\af^*_1=\af_1+x \af_6$ & 
$\af^*_2=\af_2+x \af_5$ & 
$\af^*_3=\af_3+x \af_7$ & 
$\af^*_4=\af_4+x \af_7$\\
$\af^*_5=  \alpha_5$ &
$\af^*_6=  \alpha_6$& 
$\af^*_7=  \alpha_7$.&
	\end{longtable} 
	
 Hence, $\phi\langle\0\rangle=\langle\0^*\rangle$, where $\0^*=\sum\limits_{i=1}^7 \af_i^*  \nb i.$
We are only interested in elements with $(\alpha_3, \af_4, \af_5, \af_6, \af_7)\neq (0,0,0,0,0).$
Then
\begin{enumerate}
    \item $\af_7\ne 0$. Then choosing $x=-\af_3\af_7^{-1}$, we obtain the family of representatives of distinct orbits $\langle \af  \nb 1+ \bt \nb 2+\gm \nb 4+\dt \nb 5+\eps \nb 6+\nb 7 \rangle$.
    \item $\af_7=0$ and $\af_6\ne 0$. Then choosing $x=-\af_1\af_6^{-1}$, we obtain the family of representatives of distinct orbits $\langle  \af \nb 2+\bt \nb 3+\gm \nb 4+ \dt \nb 5+ \nb 6  \rangle$.
    \item $\af_7=\af_6=0$ and $\af_5\ne 0$. Then choosing $x=-\af_2\af_5^{-1}$, we obtain the family of representatives of distinct orbits $\langle  \af \nb 1+\bt \nb 3+\gm \nb 4+  \nb 5   \rangle$.
    \item $\af_7=\af_6=\af_5=0$ and $\af_4\ne 0$. Then we obtain the family of representatives of distinct orbits $\langle  \af \nb 1+\bt \nb 2+\gm\nb 3+ \nb 4    \rangle$.
    \item $\af_7=\af_6=\af_5=\af_4=0$ and $\af_3\ne 0$. Then we obtain the family of representatives of distinct orbits $\langle  \af \nb 1+\bt \nb 2+ \nb 3     \rangle$.
\end{enumerate}
    Summarizing, we have the following distinct orbits:
\begin{longtable}{ll}

$\langle  \af \nb 1+\bt \nb 2+ \nb 3     \rangle$ &

$\langle  \af \nb 1+\bt \nb 2+\gm\nb 3+ \nb 4    \rangle$\\

$\langle  \af \nb 1+\bt \nb 3+\gm \nb 4+  \nb 5   \rangle$ &

$\langle  \af \nb 2+\bt \nb 3+\gm \nb 4+ \dt \nb 5+ \nb 6  \rangle$\\

\multicolumn{2}{c}{$\langle \af  \nb 1+ \bt \nb 2+\gm \nb 4+\dt \nb 5+\eps \nb 6+\nb 7 \rangle$}\\

\end{longtable}

This gives the following new algebras:

{\tiny \begin{longtable}{lllllllllllllll}

$\n^4_{21}(\af,\bt)$ &$:$& $e_1 e_1 = e_2$&$e_1e_2=\af e_4$&$e_1e_3=e_4$&$e_2e_1=e_3+\bt e_4$&  $e_2 e_2=e_3$ \\  \hline

$\n^4_{22}(\af,\bt,\gm)$ &$:$& $e_1 e_1 = e_2$&$e_1e_2=\af e_4$&$e_1e_3=\gm e_4$&$e_2e_1=e_3+\bt e_4$&  $e_2 e_2=e_3$&$e_3e_1=e_4$ \\  \hline

$\n^4_{23}(\af,\bt,\gm)$ &$:$& 
$e_1 e_1 = e_2$&$e_1e_2=\af e_4$&$e_1e_3=\bt e_4$&$e_2e_1=e_3$\\
&&  $e_2 e_2=e_3$&$e_2e_3=e_4$&$e_3e_1=\gm e_4$ \\  \hline

$\n^4_{24}(\af,\bt,\gm,\delta)$ &$:$& 
$e_1 e_1 = e_2$&$e_1e_3=\bt e_4$&$e_2e_1=e_3+\af e_4$&  $e_2 e_2=e_3$ \\
&&$e_2e_3=\delta e_4$&$e_3e_1=\gm e_4$ & $e_3e_2=e_4$\\  \hline

$\n^4_{25}(\af, \bt, \gm, \delta, \eps)$ &$:$& 
$e_1 e_1 = e_2$&$e_1e_2= \af e_4$& $e_2e_1=e_3+\bt e_4$ & $e_2 e_2=e_3$\\
&&$e_2e_3=\delta e_4$ & $e_3e_1=\gm e_4$ & $e_3e_2=\eps e_4$ &$e_3e_3=e_4$\\  \hline


\end{longtable}
}
	\subsection{$1$-dimensional central extensions of $\cd 3{03}$}
Here we collect all information about ${\mathfrak{CD}_{03}^{3}}$:
\begin{longtable}{|l|l|l|}
\hline
Algebra  & Multiplication & Cohomology   \\
\hline
${\mathfrak{CD}}_{03}^{3}$ &  
$\begin{array}{l}
e_1 e_1 = e_2\\ 
e_2 e_1=e_3
\end{array}$& 
$\begin{array}{lcl}
{\rm H^2_{\mathfrak{CD}}}({\mathfrak{CD}}_{03}^{3})&=&\langle  \Dl 12, \Dl 22, \Dl 13-2\Dl 31 \rangle \\
{\rm H^2}({\mathfrak{CD}}_{03}^{3})&=&{\rm H^2_{\mathfrak{CD}}}({\mathfrak{CD}_{03}^3}) \oplus 
\langle \Dl 31, \Dl 23, \Dl 32, \Dl 33    \rangle 
\end{array}$\\
\hline
\end{longtable}			
	Let us use the following notations 
	\begin{longtable}{llll}
        $\nb 1 = \Dl 12$&$ \nb 2 = \Dl 22$&  $\nb 3 = \Dl 13-2\Dl 31$ &$\nb 4=\Dl 31$\\
        $\nb 5 = \Dl 23$&$ \nb 6 = \Dl 32$&  $\nb 7 = \Dl 33$.
	\end{longtable}
	Take $\0=\sum_{i=1}^7\af_i\nb i\in {\rm H}^2(\cd 3{03})$.
	If 
	$$
	\phi=
	\begin{pmatrix}
	x &    0  &  0\\
	y &  x^2  &  0\\
	z &   xy  &  x^3
	\end{pmatrix}\in\aut{\cd 3{03}},
	$$
	then
	$$
	\phi^T\begin{pmatrix}
	0             &  \af_1  & \af_3\\
	0             &  \af_2  & \af_5\\
	\af_4-2\af_3  &    \af_6  & \af_7
	\end{pmatrix} \phi=
	\begin{pmatrix}
	\af^*                   &  \af^*_1  & \af^*_3\\
	\af^{**}                &  \af^*_2  & \af^*_5\\
	\af^*_4-2\af^*_3        &  \af^*_6  & \af^*_7
	\end{pmatrix},
	$$
	where
	\begin{longtable}{lclc}
	    $\af^*_1$ &$=$&$ (\af_1x^2 + (\af_2+\af_3)xy + \af_5y^2 + \af_6xz + \af_7yz)x,$\\
        $\af^*_2$ &$=$&$ (\af_2x^2 + (\af_5+\af_6)xy + \af_7y^2)x^2,$\\
        $\af^*_3$ &$=$&$ (\af_3x + \af_5y + \af_7z)x^3,$\\
        $\af^*_4$ &$=$&$ (\af_4x + (2\af_5+\af_6)y + 3\af_7z)x^3,$\\
        $\af^*_5$ &$=$&$ (\af_5x + \af_7y)x^4,$\\
        $\af^*_6$ &$=$&$ (\af_6x + \af_7y)x^4,$\\
        $\af^*_7$ &$=$&$ \af_7x^6.$
	\end{longtable}
	Hence, $\phi\langle\0\rangle=\langle\0^*\rangle$, where $\0^*=\sum\limits_{i=1}^7 \af_i^*  \nb i.$
We are only interested in $\theta$ with $(\af_4, \af_5, \af_6, \af_7)\neq (0,0,0,0)$.

\begin{enumerate}
    \item $\af_7\ne 0$. Then put $y=-\frac{\af_6x}{\af_7}$ and $z=\frac x{\af_7^2}(\af_5\af_6 - \af_3\af_7)$ to make $\af^*_3=\af^*_6=0$.
    \begin{enumerate}
        \item $\af_5\ne \af_6$. Then choosing $x=\frac {\af_5 - \af_6}{\af_7}$, we obtain the family of representatives of distinct orbits $\la\af\nb 1+\bt\nb 2+\gm\nb 4+\nb 5+\nb 7\ra$.
        
        \item $\af_5=\af_6$ and $\af_4\ne 3\af_3$. Then choosing $x=\sqrt{\frac {\af_4-3\af_3}{3\af_7}}$, we obtain the family of representatives $\la\af\nb 1+\bt\nb 2+\nb 4+\nb 7\ra$, where two distinct pairs $(\af,\bt)$ and $(\af',\bt')$ determine the same orbit if and only if $(\af,\bt)=(-\af',\bt')$.
        
         \item $\af_5=\af_6$, $\af_4=3\af_3$ and $\af_6^2\ne \af_2\af_7$. Then choosing $x=\frac 1{\af_7}\sqrt{\af_2\af_7-\af_6^2}$, we obtain the family of representatives $\la\af\nb 1+\nb 2+\nb 7\ra$, where two distinct parameters $\af$ and $\af'$ determine the same orbit if and only if $\af=-\af'$.
         
         \item $\af_5=\af_6$, $\af_4=3\af_3$ and $\af_6^2=\af_2\af_7$. Then we obtain 2 representatives $\la\nb 7\ra$ and $\la\nb 1+\nb 7\ra$ depending on whether $\af_3\af_6=\af_1\af_7$ or not.
    \end{enumerate} 
    
    \item $\af_7=0$, $\af_5\ne 0$ and $\af_6\ne 0$. Then $\af^*_7=0$, and we put $y=-\frac{\af_3x}{\af_5}$ and $z=\frac x{\af_5\af_6}(\af_2\af_3 - \af_1\af_5)$ to make $\af^*_1=\af^*_3=0$.
    \begin{enumerate}
        \item $(\af_4-2\af_3)\af_5 \ne \af_3\af_6$. Then choosing $x=\frac 1{\af_5\af_6}((\af_4-2\af_3)\af_5 - \af_3\af_6)$, we obtain the family of representatives of distinct orbits $\la\af\nb 2+\nb 4+\bt\nb 5+\nb 6\ra_{\bt\ne 0}$.
        \item $(\af_4-2\af_3)\af_5 = \af_3\af_6$ and $(\af_2 - \af_3)\af_5 \ne \af_3\af_6$. Then choosing $x=\frac 1{\af_5\af_6}((\af_2 - \af_3)\af_5 - \af_3\af_6)$, we obtain the family of representatives of distinct orbits $\la\nb 2+\af\nb 5+\nb 6\ra_{\af\ne 0}$.
        \item $(\af_4-2\af_3)\af_5 = \af_3\af_6=(\af_2 - \af_3)\af_5$. Then we obtain the family of representatives of distinct orbits $\la\af\nb 5+\nb 6\ra_{\af\ne 0}$.
    \end{enumerate}
    
    \item $\af_7=\af_6=0$ and $\af_5\ne 0$. Then $\af^*_6=\af^*_7=0$, and we put $y=-\frac{\af_3x}{\af_5}$ to make $\af^*_3=0$.
    \begin{enumerate}
        \item $\af_4\ne 2\af_3$. Then choosing $x=\frac {\af_4-2\af_3}{\af_5}$, we obtain the family of representatives of distinct orbits $\la\af\nb 1+\bt\nb 2+\nb 4+\nb 5\ra$.
        \item $\af_4=2\af_3$ and $\af_2\ne\af_3$. Then choosing $x=\frac {\af_2-\af_3}{\af_5}$, we obtain the family of representatives of distinct orbits $\la\af\nb 1+\nb 2+\nb 5\ra$.
        \item $\af_4=2\af_3$ and $\af_2=\af_3$. Then we obtain 2 representatives $\la\nb 5\ra$ and $\la\nb 1+\nb 5\ra$ depending on whether $\af_3^2 = \af_1\af_5$ or not.
    \end{enumerate}
    
    \item $\af_7=\af_5=0$ and $\af_6\ne 0$. Then $\af^*_5=\af^*_7=0$, and we put $y=-\frac{\af_4x}{\af_5}$ and $z=\frac x{\af_6^2}((\af_2 + \af_3)\af_4 - \af_1\af_6)$ to make $\af^*_1=\af^*_4=0$.
    \begin{enumerate}
        \item $\af_2\ne\af_4$. Then choosing $x=\frac {\af_2-\af_4}{\af_6}$, we obtain the family of representatives of distinct orbits $\la\nb 2+\af\nb 3+\nb 6\ra$. 
        \item $\af_2=\af_4$. Then we obtain 2 representatives $\la\nb 6\ra$ and $\la\nb 3+\nb 6\ra$ depending on whether $\af_3=0$ or not. We will join $\la\nb 6\ra$ with the family $\la\af\nb 5+\nb 6\ra_{\af\ne 0}$ found above.
    \end{enumerate}
    
     \item $\af_7=\af_6=\af_5=0$ and $\af_4\ne 0$. Then $\af^*_5=\af^*_6=\af^*_7=0$.
     \begin{enumerate}
         \item $\af_2\ne-\af_3$. Then choosing $y = -\frac{\af_1x}{\af_2 + \af_3}$, we obtain the family of representatives of distinct orbits $\la\af\nb 2+\bt\nb 3+\nb 4\ra_{\bt\ne-\af}$. 
         \item $\af_2=-\af_3$. Then we obtain two families of representatives of distinct orbits $\la\af\nb 2-\af\nb 3+\nb 4\ra$ and $\la\nb 1+\af\nb 2-\af\nb 3+\nb 4\ra$ depending on whether $\af_1=0$ or not. The family $\la\af\nb 2-\af\nb 3+\nb 4\ra$ will be joined with the family $\la\af\nb 2+\bt\nb 3+\nb 4\ra_{\bt\ne-\af}$ from the previous item.
         
     \end{enumerate}
    \end{enumerate}
    
    Summarizing, we have the following distinct orbits:
\begin{longtable}{lll}
    $\la\nb 1+\af\nb 2-\af\nb 3+\nb 4\ra$ &
    $\la\af\nb 1+\bt\nb 2+\nb 4+\nb 5\ra$ &
    $\la\af\nb 1+\bt\nb 2+\gm\nb 4+\nb 5+\nb 7\ra$ \\
    \multicolumn{2}{l}{$\la\af\nb 1+\bt\nb 2+\nb 4+\nb 7\ra^{O(\af,\bt)=O(-\af,\bt)}$}&
    $\la\af\nb 1+\nb 2+\nb 5\ra$ \\
    $\la\af\nb 1+\nb 2+\nb 7\ra^{O(\af)=O(-\af)}$ &
    $\la\nb 1+\nb 5\ra$ &
    $\la\nb 1+\nb 7\ra$ \\
    $\la\af\nb 2+\bt\nb 3+\nb 4\ra$ &  
    $\la\nb 2+\af\nb 3+\nb 6\ra$&
    $\la\af\nb 2+\nb 4+\bt\nb 5+\nb 6\ra_{\bt\ne 0}$ \\
    $\la\nb 2+\af\nb 5+\nb 6\ra_{\af\ne 0}$ &
    $\la\nb 3+\nb 6\ra$ &
    $\la\nb 5\ra$ \\
    $\la\af\nb 5+\nb 6\ra$&
    $\la\nb 7\ra$
\end{longtable}

They correspond to the following new algebras:
{\tiny 
\begin{longtable}{lllllllllllllll}

$\n^4_{26}(\af)$ &$:$& 
$e_1 e_1 = e_2$&$e_1e_2=e_4$&$e_1e_3=-\af e_4$ &$e_2e_1=e_3$  &$e_2e_2=\af e_4$&$ e_3e_1=(1+2 \af) e_4$  \\  \hline
 
$\n^4_{27}(\af,\bt)$ &$:$& $e_1 e_1 = e_2$&$e_1e_2=\af e_4$&$e_2e_1=e_3$ &$e_2e_2=\bt e_4$&$e_2e_3=e_4$&$e_3e_1=e_4$&  \\  \hline

$\n^4_{28}(\af,\bt,\gm)$ &$:$& $e_1 e_1 = e_2$&$e_1e_2=\af e_4$&$e_2e_1=e_3$&$e_2e_2=\bt e_4$\\&&$e_2e_3=e_4$&$e_3e_1=\gm e_4$&$e_3e_3=e_4$&  \\  \hline

$\n^4_{29}(\af,\bt)$ &$:$& $e_1 e_1 = e_2$&$e_1e_2=\af e_4$&$e_2e_2=\bt e_4$ &$e_2e_1=e_3$&$e_3e_1=e_4$&$e_3e_3=e_4$&  \\  \hline

$\n^4_{30}(\af) $ &$:$& $e_1 e_1 = e_2$&$e_1e_2=\af e_4$&$e_2e_1=e_3$&$e_2e_2=e_4$&$e_2e_3=e_4$  \\  \hline

$\n^4_{31}(\af) $ &$:$& $e_1 e_1 = e_2$&$e_1e_2=\af e_4$&$e_2e_1=e_3$&$e_2e_2=e_4$&$e_3e_3=e_4$  \\  \hline

$\n^4_{32}$ &$:$& $e_1 e_1 = e_2$&$e_1e_2=e_4$&$e_2e_1=e_3$&$e_2e_3=e_4$  \\  \hline

$\n^4_{33}$ &$:$& $e_1 e_1 = e_2$&$e_1e_2=e_4$&$e_2e_1=e_3$&$e_3e_3=e_4$  \\  \hline

$\n^4_{34}(\af,\bt)$ &$:$& $e_1 e_1 = e_2$&$e_1e_3=\bt e_4$&$e_2e_1=e_3$&$e_2e_2=\af e_4$&
\multicolumn{2}{l}{$e_3e_1=(1-2\bt) e_4$}  \\  \hline

$\n^4_{35}(\af)$ &$:$& $e_1 e_1 = e_2$&$e_1e_3=\af e_4$&$e_2e_1=e_3$ &$e_2e_2=e_4$&$e_3e_1=-2\af e_4$&$e_3e_2=e_4$&  \\  \hline

$\n^4_{36}(\af,\bt)_{\bt\neq 0}$ &$:$& $e_1 e_1 = e_2$&$e_2e_1=e_3$&$e_2e_2=\af e_4$ 
&$e_2e_3=\bt e_4$&$e_3e_1=e_4$&$e_3e_2=e_4$&  \\  \hline

$\n^4_{37}(\af)_{\af\neq0}$ &$:$& $e_1 e_1 = e_2$&$e_2e_1=e_3$&$e_2e_2=e_4$&$e_2e_3=\af e_4$&$e_3e_2=e_4$  \\  \hline

$\n^4_{38}$ &$:$& $e_1 e_1 = e_2$&$e_1e_3=e_4$&$e_2e_1=e_3$&$e_3e_1=-2e_4$&$e_3e_2=e_4$ \\  \hline

$\n^4_{39}$ &$:$& $e_1 e_1 = e_2$&$e_2e_1=e_3$&$e_2e_3=e_4$&  \\  \hline

$\n^4_{40}(\af)$ &$:$& $e_1 e_1 = e_2$&$e_2e_1=e_3$&$e_2e_3=\af e_4$&$e_3e_2=e_4$&  \\  \hline

$\n^4_{41}$ &$:$& $e_1 e_1 = e_2$&$e_2e_1=e_3$&$e_3e_3=e_4$&  \\  \hline


\end{longtable}
 }

 \begin{longtable}{ll}
$\n^4_{29}(\af,\bt) \cong \n^4_{29}(-\af,\bt)$ &
$\n^4_{31}(\af) \cong \n^4_{31}(-\af)$
 \end{longtable}
 
  \subsection{$1$-dimensional central extensions of $\cd 3{04}$}\label{sec-CD^3_04}
Here we  collect all information about ${\mathfrak{CD}_{04}^{3}}:$
\begin{longtable}{|l|l|l|}
\hline
Algebra  & Multiplication & Cohomology   \\
\hline
${\mathfrak{CD}}_{04}^{3}$ &  
$\begin{array}{l}
 e_1 e_1 = e_2\\
 e_1 e_2=e_3\\ 
 e_2 e_1=\lambda e_3
\end{array}$& 
$\begin{array}{lcl}
{\rm H^2_{\mathfrak{CD}}}({\mathfrak{CD}}_{04}^{3})&=&\langle  
 (\lb-2)\Dt 13-(2\lb-1)\Dt 31,  \Dt 21, \Dt 22
     \rangle \\
{\rm H^2}({\mathfrak{CD}}_{04}^{3})&=&{\rm H^2_{\mathfrak{CD}}}({\mathfrak{CD}_{04}^3}) \oplus 
\langle
\Dl 13-2 \Dl 31,  \Dl 23,  \Dl 32, \Dl 33
\rangle 
\end{array}$\\
\hline
\end{longtable}	
Let us use the following notations 
\begin{longtable}{llll}
$\nb 1 = (\lb-2)\Dl 13-(2\lb-1)\Dl 31$ &$\nb 2 = \Dl 21$ & $\nb 3 = \Dl 22$\\
$\nb 4= \Dl 13-2 \Dl 31$ & $\nabla_5=\Dl 23$ & $\nabla_6=\Dl 32$ & $\nabla_7=\Dl 33$
\end{longtable}

Take $\0=\sum_{i=1}^3\af_i\nb i\in {\rm H}^2_{\mathfrak{CD}}(\cd 3{04})$.
If 
$$
\phi=
\begin{pmatrix}
x &               0  &  0\\
y &             x^2  &  0\\
z &   (\lambda+1)xy  &  x^3
\end{pmatrix}\in\aut{\cd 3{04}},
$$
then
$$
\phi^T\begin{pmatrix}
0               &      0  & (\lb-2)\af_1+\af_4\\
\af_2           &  \af_3  & \af_5\\
-(2\lb-1)\af_1-2 \af_4  &     \af_6   & \af_7
\end{pmatrix} \phi=
	\begin{pmatrix}
	\af^*                &  \af^{**} & (\lb-2)\af_1^*+\af^*_4\\
	\af_2^*+\lb\af^{**}  &  \af_3^*  & \af^*_5\\
	-(2\lb-1)\af_1^*-2 \af^*_4     &  \af^*_6  & \af^*_7
\end{pmatrix},
$$
where
\begin{longtable}{lll}
    $\af^*_1$&  $=$& $\frac{x^3}{3}  (3 x \af_1-y (2 \af_5+\af_6)-3 z \af_7)$\\ 
    $\af^*_2$& $=$& $x (x^2 \af_2+y (1+\lb ) (z \af_7 (1-\lb )+y (\af_6-\af_5 \lb ))+x (z (\af_5- \af_6 \lb) -y (\af_3 (\lb-1 )+$\\
    &&$\af_1 (\lb-1 ) (\lb+1 )^2+\af_4 (2+3 \lb +\lb ^2))))$\\
    $\af_3^*$ &$=$&$x^2 (x^2 \af_3+x y (\af_5+\af_6) (1+\lb )+y^2 \af_7 (\lb+1 )^2)$\\
    $\af_4^*$ &$=$&$\frac{ x^3}{3} (3 x \af_4+y \af_6 (\lb-2 )+3 z \af_7 (\lb-1 )+y \af_5 (2 \lb-1 ))$\\
    $\af_5^*$ &$=$&$x^4 (x \af_5+y \af_7 (\lb+1 ))$\\
    $\af_6^*$ &$=$&$x^4 (x \af_6+y \af_7 (\lb+1 ))$\\
    $\af_7^*$ &$=$&$x^6 \af_7$
\end{longtable} 
Hence, $\phi\langle\0\rangle=\langle\0^*\rangle$, where $\0^*=\sum\limits_{i=1}^7 \af_i^*  \nb i.$ We are only interested in elements with $(\alpha_4, \af_5, \af_6,\af_7)\neq (0,0,0,0).$ 

\begin{enumerate}
    \item $\alpha_7\ne 0$ and $\lambda\neq -1$. We may assume $\alpha_7=1$. We put $z = \af_1x - \frac 13(2\af_5 + \af_6)y$ to make $\af^*_1=0$.
\begin{enumerate}
    \item 
    $\af_6\neq \af_5$. 
    Then choosing $y=-\frac{3((\lb - 1)\af_1 + \af_4)x}{\af_5 - \af_6}$ we have $\af^*_4=0$. Now, if $\af_6\ne\frac{3(\lb+1)((\lb-1)\af_1 + \af_4)}{\af_5 - \af_6}$, then choosing $x=\af_6-\frac{3(\lb+1)((\lb-1)\af_1 + \af_4)}{\af_5 - \af_6}$
    we have the family of representatives of distinct orbits
    $\la \af \nb 2+\bt \nb 3+\gm \nb 5+ \nb 6+ \nb 7\ra_{\lb\ne-1,\gm\ne 1}$. Otherwise, choosing $x=\af_5-\af_6$ we have the family of representatives of distinct orbits
    $\la \af \nb 2+\bt \nb 3+\nb 5+ \nb 7\ra_{\lb\ne-1}$.
    
    \item 
   $\af_6= \af_5$ and $\af_4 \neq (1 - \lb )\af_1.$ 
    Then choosing $y=-\frac{\af_5x}{\lb + 1}$ and
    $x=\sqrt{\af_4+(\lb-1)\af_1}$,
    we have the family of representatives 
    $\la \af \nb 2+ \bt \nb 3+  \nb 4+ \nb 7\ra_{\lb\ne-1}$, where $(\af,\bt)$ and $(\af',\bt')$ determine the same orbit if and only if $(\af,\bt)=(\pm\af',\bt')$.

    \item 
    $\af_6= \af_5,$ $\af_4 = (1 - \lb )\af_1$ and  $\af_3\ne \af_5^2.$ 
    Then choosing  $y=-\frac{\af_5x}{\lb + 1}$ and
        $x=\sqrt{ \af_3-\af_5^2}$,
    we have the family of representatives 
    $\la \af \nb 2+  \nb 3+ \nb 7\ra_{\lb\ne-1}$, where $\af$ and $\af'$ determine the same orbit if and only if $\af=\pm\af'$.

  \item 
      $\af_6= \af_5,$ $\af_4=(1-\lb)\af_1$ and  $\af_3=\af_5^2$.  
    Then have two   representatives 
    $\la  \nb 7\ra_{\lb\ne-1}$  and  $\la \nb 2+ \nb 7\ra_{\lb\ne-1}$ depending on whether $\af_2=(\lb -1)\af_1 \af_5$ or not.
\end{enumerate}    
\item $\af_7\ne 0$ and $\lambda=-1$. We may assume $\alpha_7=1$.
\begin{enumerate}
    \item $\af_6\neq0$ and $\af_6\neq \af_5.$
        Then choosing 
    $x=\af_6,$
    $y=\frac{3(2 \af_1-\af_4) x}{\af_5-\af_6}$
    and 
    $z=\frac{(\af_4 (2 \af_5+\af_6)-3 \af_1 (\af_5+\af_6))x}{\af_5-\af_6},$
    we have the family of representatives of distinct orbits
    $\la \af \nb 2+\bt \nb 3+\gm \nb 5+\nb 6+ \nb 7\ra_{\lb=-1, \gm \neq 1}$.

     \item $\af_6\neq0, \af_6= \af_5$ and $\af_3\neq \af_5^2$.
        Then choosing 
    $x=\af_5,$
    $y=\frac{(\af_2+2 \af_1 \af_5)x}{2 (\af_5^2-\af_3)}$
    and 
    $z=\frac{(\af_2 \af_5+2 \af_1 \af_3)x}{2 (\af_3-\af_5^2)},$
    we have the family of representatives of distinct orbits
    $\la  \af \nb 3+ \bt \nb 4+ \nb 5+\nb 6+ \nb 7\ra_{\lb=-1,\af\ne 1}$.

     \item $\af_6\neq0, \af_6= \af_5$ and $\af_3= \af_5^2.$
        Then choosing 
    $x=\af_6,$
    $y=\af_1$
    and 
    $z=0,$
    we have the family of representatives of distinct orbits
    $\la \af \nb 2+ \nb 3+ \bt \nb 4+ \nb 5+\nb 6+ \nb 7\ra_{\lb=-1}$.
   
     \item $\af_6=0$ and $\af_5\neq 0.$
        Then choosing 
    $x=\af_5,$
    $y=3(2\af_1-\af_4)$
    and 
    $z=(2 \af_4-3 \af_1)\af_5,$
    we have the family of representatives of distinct orbits
    $\la \af \nb 2+ \bt \nb 3+  \nb 5 + \nb 7\ra_{\lb=-1}$.

   \item $\af_6=0,$  $\af_5= 0,$ $\af_4\neq 2 \af_1$ and $\af_3\neq 0.$
        Then choosing 
    $x=\sqrt{\af_4-2 \af_1},$
    $y=-\frac{\af_2x}{2\af_3}$
    and 
    $z=\af_1x,$
    we have the family of representatives of distinct orbits
    $\la \af  \nb 3+  \nb 4+   \nb 7\ra_{\lb=-1, \af\ne 0}$.

  \item $\af_6=0,$  $\af_5= 0,$ $\af_4= 2 \af_1$ and $\af_3\neq 0.$
        Then choosing 
    $x=\sqrt{\af_3},$
    $y=-\frac{\af_2x}{2 \af_3}$
    and 
    $z=\af_1x,$
    we have the representative
    $\la  \nb 3+\nb 7\ra_{\lb=-1}$.

  \item $\af_6=0,$  $\af_5= 0,$ $\af_4= 2 \af_1$ and $\af_3= 0.$
        Then choosing 
    $z=\af_1x$,
    we have two representatives 
    $\la \nb 7\ra_{\lb=-1}$ and $\la  \nb 2+ \nb 7\ra_{\lb=-1}$, 
    depending on whether $\af_2=0$ or not.

\end{enumerate}

\item $\af_7=0$ and $\af_6\ne 0.$ We may assume $\alpha_6=1$.
\begin{enumerate}
    \item $ \af_5\not\in\{-\frac 12,\lb\}$. Then the system $\af^*_1=\af^*_2=0$ has a unique solution in $y$ and $z$. Now, if $\af_4\neq -\frac{((2 \lb-1)\af_5+\lb-2)\af_1}{2\af_5 + 1}$, then choosing the suitable value of $x$
    we have the family of representatives distinct orbits
     $\la \af \nb 3+ \nb 4 +\bt \nabla_5+ \nb 6\ra_{\bt\not\in\{-\frac 12,\lb\}}$. Otherwise, we have two families of representatives distinct orbits $\la \af \nabla_5+ \nb 6\ra_{\af\not\in\{-\frac 12,\lb\}}$ and $\la \nb 3+ \af \nabla_5+\nb 6\ra_{\af\not\in\{-\frac 12,\lb\}}$ depending on whether $\af_3=-\frac{3(\lb + 1)(\af_5 + 1)\af_1}{2\af_5 + 1}$ or not.

     
       \item $\af_5= \lb\ne-\frac 12$ and $\af_4  \neq \frac{2 (1-\lb ^2)\af_1}{2\lb+1}$. 
    Then choosing
    $x=\af_4+\frac{2 (\lb^2-1)\af_1}{2 \lb+1}$ and
    $y=\frac{3 \af_1 x}{2 \lb+1}$  
    we have the family of representatives of distinct orbits
     $\la  \af \nb 2+\bt \nb 3+ \nb 4+\lb \nabla_5+\nb 6\ra_{\lb\ne -\frac 12}$.

    \item $\af_5= \lb\ne-\frac 12$ and  $\af_4 = \frac{2 (1-\lb ^2)\af_1}{2\lb+1}$. We put $y=\frac{3 \af_1 x}{2 \lb+1}$ to make $\af^*_1=0$. Now, if $\af_3\neq -\frac{3 (\lb+1)^2\af_1}{2 \lb + 1}$,
    then choosing
    $x=\af_3+\frac{3(\lb+1)^2\af_1}{2 \lb+1}$, 
    we have the family of representatives of distinct orbits
     $\la  \af \nb 2+  \nb 3 +\lb \nabla_5+ \nb 6\ra_{\lb\ne -\frac 12}$. Otherwise,  
    we have two representatives
     $\la \lb \nabla_5+\nb 6\ra_{\lb\ne -\frac 12}$ and $\la   \nb 2+  \lb \nabla_5+\nb 6\ra_{\lb\ne -\frac 12}$,
    depending on whether $\af_2=-\frac{9(\lb - 1)(\lb + 1)^2\af_1^2}{(2\lb + 1)^2}$ or not.
     
    \item $\af_5=-\frac 12 \ne\lb$.
    Then put
    $y=2 \af_4 x$
    and 
    $z=\frac{2(\af_2-2 \af_4 (\lb-1) (\af_3+\af_1 (\lb+1)^2)x}{2 \lb+1}$
    to make $\af^*_2=\af^*_4=0$. Now, if $\af_3\ne -(\lb+1)\af_4$, then choosing $x=\af_3+(\lb+1)\af_4,$
    we have the family of representatives  of distinct orbits
        $\la \af \nb 1+ \nb 3 -\frac 12\nb 5+ \nabla_6\ra_{\lb\ne -\frac 12}$. Otherwise, we have two  representatives
         $\la -\frac 12\nb 5+  \nabla_6\ra_{\lb\ne -\frac 12}$ and $\la \nb 1 -\frac 12\nb 5+  \nabla_6\ra_{\lb\ne -\frac 12}$,
         depending on whether $\af_1=0$ or not. 
        

 \item $\af_5=-\frac 12=\lb$ and $\af_4\ne-2 \af_3$.
    Then choosing
    $x=\af_3+\frac{\af_4}2$ and
    $y=2\af_4x$
    we have the family of representatives of distinct orbits
         $\la \af\nabla_1 +\beta \nabla_2+\nabla_3 -\frac 12\nb 5+  \nabla_6\ra_{\lb= -\frac 12}$.

 \item  $\af_5=-\frac 12=\lb$ and $\af_4=-2 \af_3$. Then we put $y=-4\af_3 x$ to make $\af^*_4=0$. Now, if $\af_2\neq\frac {3\af_3}2 (\af_1+4 \af_3)$,
    then choosing
    $x=\sqrt{\af_2-\frac{3\af_3}2 (\af_1+4 \af_3)}$,
    we have the family of representatives of distinct orbits
         $\la \af\nabla_1 +  \nabla_2 -\frac 12\nb 5+  \nabla_6\ra_{\lb= -\frac 12}$. Otherwise, we have two  representatives
    $\la \nabla_1 -\frac 12\nb 5+  \nabla_6\ra_{\lb= -\frac 12}$ and 
    $\la  -\frac 12\nb 5+  \nabla_6\ra_{\lb= -\frac 12}$,
    depending on whether $\af_1=0$ or not. 


\end{enumerate}
    
\item $\af_7=\af_6=0$ and $\af_5\ne 0$. We may assume $\af_5=1$. Then the system $\af^*_1=\af^*_2=0$ has a unique solution in $y$ and $z$. Now, if $2 \af_4+(2 \lb-1)\af_1\neq 0$, then choosing $x=\af_4+(\lb-\frac 12)\af_1$ we have the family of representatives of distinct orbits $\la \af \nabla_3+ \nabla_4+\nb 5\ra$. Otherwise, we have two representatives $\la  \nb 5\ra$ and $\la \nabla_3+ \nb 5\ra$ depending on whether $2 \af_3+3(\lb+1) \af_1=0$ or not.

\item $\af_7=\af_6=\af_5=0$ and $\af_4\ne 0$. We may assume $\af_4=1$.

\begin{enumerate}
    \item $(\lb-1) ((\lb+1)^2\af_1+\af_3)+(\lb + 1)(\lb + 2)\neq 0.$
    Then choosing 
    $y=\frac{\af_2x}{(\lb-1) ((\lb+1)^2\af_1+\af_3)+(\lb + 1)(\lb + 2)}$
    we have the family of representatives of distinct orbits
         $\la \af \nabla_1+ \bt \nabla_3+\nb 4\ra$, where $(\lb-1) ((\lb+1)^2\af+\bt)+(\lb + 1)(\lb + 2)\neq 0$.

 \item $(\lb-1) ((\lb+1)^2\af_1+\af_3)+(\lb + 1)(\lb + 2)= 0$. Then $\lb\neq 1$
    and we have two families of representatives of distinct orbits $\La \af \nabla_1+ \frac{(\lb+1) ((\lb^2-1)\af + \lb + 2)}{1-\lb}  \nabla_3+\nb 4\Ra_{\lb\ne 1}$ and
         $\La \af \nabla_1+\nabla_2+ \frac{(\lb+1) ((\lb^2-1)\af + \lb + 2)}{1-\lb}  \nabla_3+\nb 4\Ra_{\lb\ne 1}$, depending on whether $\af_2=0$ or not. The first family will be joined with the family from the previous item.
         
\end{enumerate}

\end{enumerate}

Summarizing, we have the following representatives of distinct orbits:

\begin{longtable}{lllllll}
$\La \af \nabla_1+\nabla_2+ \frac{(\lb+1) ((\lb^2-1)\af + \lb + 2)}{1-\lb}  \nabla_3+\nb 4\Ra_{\lb\ne 1}$ &
$\la \af\nabla_1 +\beta \nabla_2+\nabla_3 -\frac 12\nb 5+  \nabla_6\ra_{\lb= -\frac 12}$\\
$\la \af\nabla_1 +  \nabla_2 -\frac 12\nb 5+  \nabla_6\ra_{\lb= -\frac 12}$&
$\la \af \nabla_1+ \bt \nabla_3+\nb 4\ra$\\
$\la \af \nb 1+ \nb 3 -\frac 12\nb 5+ \nabla_6\ra_{\lb\ne -\frac 12}$&
$\la \nabla_1 -\frac 12\nb 5+  \nabla_6\ra_{\lb= -\frac 12}$\\
$\la  \af \nb 2+\bt \nb 3+ \nb 4+\lb \nabla_5+\nb 6\ra_{\lb\ne -\frac 12}$&
$\la \af \nb 2+ \nb 3+ \bt \nb 4+ \nb 5+\nb 6+ \nb 7\ra_{\lb=-1}$\\
$\la \af \nb 2+ \bt \nb 3+ \nb 4+ \nb 7\ra_{\lb\ne-1}^{O(\af,\bt)=O(-\af,\bt)}$&
$\la  \af \nb 2+  \nb 3 +\lb \nabla_5+ \nb 6\ra_{\lb\ne -\frac 12}$\\
$\la \af \nb 2+\bt \nb 3+\gm \nb 5+ \nb 6+ \nb 7\ra_{\gm\ne 1}$&
$\la \af \nb 2+\bt \nb 3+\nb 5+ \nb 7\ra$\\
$\la \af \nb 2+  \nb 3+ \nb 7\ra_{\lb\ne-1}^{O(\af)=O(-\af)}$&
$\la   \nb 2+  \lb \nabla_5+\nb 6\ra_{\lb\ne -\frac 12}$\\
$\la \nb 2+ \nb 7\ra$&
$\la \af \nabla_3+ \nabla_4+\nb 5\ra$\\
$\la \af \nb 3+ \nb 4 +\bt \nabla_5+ \nb 6\ra_{\bt\not\in\{-\frac 12,\lb\}}$&
$\la \af \nb 3+ \bt \nb 4+ \nb 5+\nb 6+ \nb 7\ra_{\lb=-1,\af\ne 1}$\\
$\la \af \nb 3+  \nb 4+   \nb 7\ra_{\lb=-1, \af\ne 0}$&
$\la \nabla_3+ \nb 5\ra$\\
$\la \nb 3+ \af \nabla_5+\nb 6\ra_{\af\not\in\{-\frac 12,\lb\}}$&
$\la  \nb 3+\nb 7\ra_{\lb=-1}$\\
$\la  \nb 5\ra$&
$\la \af \nabla_5+ \nb 6\ra$\\
$\la  \nb 7\ra$ &

\end{longtable}
 
The corresponding algebras are:
{\tiny \begin{longtable}{lllllllllllllll}

$\n^4_{42}(\lb,\af)_{\lb\ne 1}$ & $:$ & $e_1 e_1 = e_2$ & $e_1 e_2=e_3$ && $e_1e_3=(\af(\lb-2)+1)e_4$\\
&& $e_2 e_1=\lb e_3+e_4$& \multicolumn{2}{l}{$e_2e_2=\frac{(\lb+1) ((\lb^2-1)\af + \lb + 2)}{1-\lb}e_4$} & $e_3e_1=(\af(1-2\lb)-2)e_4$  \\ 

\hline

$\n^4_{43}(\af,\bt)$ & $:$ & $e_1 e_1 = e_2$ & $e_1 e_2=e_3$ & $e_1e_3=-\frac {5\af}2e_4$ & $e_2 e_1=-\frac 12e_3+\bt e_4$\\
&& $e_2e_2=e_4$ & $e_2 e_3=-\frac 12 e_4$ & $e_3e_1=2\af e_4$ & $e_3e_2 = e_4$ \\

\hline

$\n^4_{44}(\af)$ & $:$ & $e_1 e_1 = e_2$ & $e_1 e_2=e_3$ & $e_1e_3=-\frac {5\af}2e_4$ & $e_2 e_1=-\frac 12e_3+e_4$\\
&& $e_2 e_3=-\frac 12 e_4$ & $e_3e_1=2\af e_4$ & $e_3e_2 = e_4$ \\

\hline

$\n^4_{45}(\lb,\af,\bt)$ & $:$ & $e_1 e_1 = e_2$ & $e_1 e_2=e_3$ & \multicolumn{2}{l}{$e_1e_3=(\af(\lb-2)+1)e_4$}\\
&& $e_2 e_1=\lb e_3$& $e_2e_2=\bt e_4$ & \multicolumn{2}{l}{$e_3e_1=(\af(1-2\lb)-2)e_4$}  \\ 

\hline

$\n^4_{46}(\lb,\af)_{\lb\ne-\frac 12}$ & $:$ & $e_1 e_1 = e_2$ & $e_1 e_2=e_3$ & $e_1e_3=\af(\lb-2)e_4$ & $e_2 e_1=\lb e_3$\\
&& $e_2e_2=e_4$ &  $e_2 e_3=-\frac 12 e_4$ & $e_3e_1=\af(1-2\lb)e_4$ & $e_3e_2 = e_4$  \\ 

\hline

$\n^4_{47}$ & $:$ & $e_1 e_1 = e_2$ & $e_1 e_2=e_3$ & $e_1e_3=-\frac 52e_4$ & $e_2 e_1=-\frac 12e_3$\\
&& $e_2 e_3=-\frac 12 e_4$ & $e_3e_1=2 e_4$ & $e_3e_2 = e_4$ \\

\hline

$\n^4_{48}(\lb,\af,\bt)_{\lb\ne -\frac 12}$ & $:$ & $e_1 e_1 = e_2$ & $e_1 e_2=e_3$ & $e_1e_3=e_4$ & $e_2 e_1=\lb e_3+\af e_4$\\
&& $e_2e_2=\bt e_4$ & $e_2 e_3=\lb e_4$ & $e_3e_1=-2e_4$ & $e_3e_2 = e_4$ \\

\hline

$\n^4_{49}(\af,\bt)$ & $:$ & $e_1 e_1 = e_2$ & $e_1 e_2=e_3$ & $e_1e_3=\bt e_4$\\
&& $e_2 e_1=-e_3+\af e_4$ & $e_2e_2=e_4$ & $e_2 e_3=e_4$\\
&& $e_3e_1=-2\bt e_4$ & $e_3e_2 = e_4$ & $e_3e_3=e_4$ \\

\hline

$\n^4_{50}(\lb,\af,\bt)_{\lb\ne-1}$ & $:$ & $e_1 e_1 = e_2$ & $e_1 e_2=e_3$ & $e_1e_3=e_4$ & $e_2 e_1=\lb e_3+\af e_4$\\
&& $e_2e_2=\bt e_4$ & $e_3e_1=-2e_4$  & $e_3e_3=e_4$ \\

\hline

$\n^4_{51}(\lb,\af)_{\lb\ne -\frac 12}$ & $:$ & $e_1 e_1 = e_2$ & $e_1 e_2=e_3$  & $e_2 e_1=\lb e_3+\af e_4$\\
&& $e_2e_2=e_4$  & $e_2e_3=\lb e_4$  & $e_3e_2=e_4$ \\

\hline

$\n^4_{52}(\lb,\af,\bt,\gm)_{\gm\ne 1}$ & $:$ & $e_1 e_1 = e_2$ & $e_1 e_2=e_3$  & $e_2 e_1=\lb e_3+\af e_4$ & $e_2e_2=\bt e_4$\\
&& $e_2e_3=\gm e_4$  & $e_3e_2=e_4$ & $e_3e_3=e_4$ \\

\hline

$\n^4_{53}(\lb,\af,\bt)$ & $:$ & $e_1 e_1 = e_2$ & $e_1 e_2=e_3$  & $e_2 e_1=\lb e_3+\af e_4$\\
&& $e_2e_2=\bt e_4$ & $e_2e_3=e_4$  & $e_3e_3=e_4$ \\

\hline

$\n^4_{54}(\lb,\af)_{\lb\ne -1}$ & $:$ & $e_1 e_1 = e_2$ & $e_1 e_2=e_3$  & $e_2 e_1=\lb e_3+\af e_4$\\
&& $e_2e_2=e_4$  & $e_3e_3=e_4$ \\

\hline

$\n^4_{55}(\lb)_{\lb\ne -\frac 12}$ & $:$ & $e_1 e_1 = e_2$ & $e_1 e_2=e_3$  & $e_2 e_1=\lb e_3+e_4$\\
&& $e_2e_3=\lb e_4$  & $e_3e_2=e_4$ \\

\hline

$\n^4_{56}(\lb)$ & $:$ & $e_1 e_1 = e_2$ & $e_1 e_2=e_3$  & $e_2 e_1=\lb e_3+e_4$ & $e_3e_3=e_4$ \\

\hline 

$\n^4_{57}(\lb,\af)$ & $:$ & $e_1 e_1 = e_2$ & $e_1 e_2=e_3$ & $e_1e_3=e_4$ & $e_2 e_1=\lb e_3$\\
&& $e_2e_2=\af e_4$ & $e_2e_3=e_4$ & $e_3e_1=-2e_4$ \\

\hline

$\n^4_{58}(\lb,\af,\bt)_{\bt\not\in\{-\frac 12,\lb\}}$ & $:$ & $e_1 e_1 = e_2$ & $e_1 e_2=e_3$ & $e_1e_3=e_4$ & $e_2 e_1=\lb e_3$\\
&& $e_2e_2=\af e_4$ & $e_2e_3=\bt e_4$ & $e_3e_1=-2e_4$ & $e_3e_2=e_4$ \\

\hline

$\n^4_{59}(\af,\bt)_{\af\ne 1}$ & $:$ & $e_1 e_1 = e_2$ & $e_1 e_2=e_3$ & $e_1e_3=\bt e_4$\\
&& $e_2 e_1=-e_3$ & $e_2e_2=\af e_4$ & $e_2e_3=e_4$\\
&& $e_3e_1=-2\bt e_4$ & $e_3e_2=e_4$ & $e_3e_3=e_4$ \\

\hline

$\n^4_{60}(\af)_{\af\ne 0}$ & $:$ & $e_1 e_1 = e_2$ & $e_1 e_2=e_3$ & $e_1e_3=e_4$ & $e_2 e_1=-e_3$\\
&& $e_2e_2=\af e_4$ & $e_3e_1=-2e_4$ & $e_3e_3=e_4$ \\

\hline

$\n^4_{61}(\lb)$ & $:$ & $e_1 e_1 = e_2$ & $e_1 e_2=e_3$ & $e_2 e_1=\lb e_3$\\
&& $e_2e_2=e_4$ & $e_2e_3=e_4$ \\

\hline

$\n^4_{62}(\lb,\af)_{\af\not\in\{-\frac 12,\lb\}}$ & $:$ & $e_1 e_1 = e_2$ & $e_1 e_2=e_3$ & $e_2 e_1=\lb e_3$\\
&& $e_2e_2=e_4$ & $e_2e_3=\af e_4$ & $e_3e_2=e_4$ \\

\hline

$\n^4_{63}$ & $:$ & $e_1 e_1 = e_2$ & $e_1 e_2=e_3$ & $e_2 e_1=-e_3$\\
&& $e_2e_2=e_4$ & $e_3e_3=e_4$ \\

\hline

$\n^4_{64}(\lb)$ & $:$ & $e_1 e_1 = e_2$ & $e_1 e_2=e_3$ & $e_2 e_1=\lb e_3$ & $e_2e_3=e_4$ \\

\hline

$\n^4_{65}(\lb,\af)$ & $:$ & $e_1 e_1 = e_2$ & $e_1 e_2=e_3$ & $e_2 e_1=\lb e_3$\\
&& $e_2e_3=\af e_4$ & $e_3e_2=e_4$ \\

\hline

$\n^4_{66}(\lb)$ & $:$ & $e_1 e_1 = e_2$ & $e_1 e_2=e_3$ & $e_2 e_1=\lb e_3$ & $e_3e_3=e_4$\\ 

\hline
\end{longtable}}\begin{longtable}{ll}
$\n^4_{50}(\lb,\af,\bt) \cong \n^4_{50}(\lb,-\af,\bt)$ &
$\n^4_{54}(\lb,\af) \cong \n^4_{54}(\lb,-\af)$
 \end{longtable}
 
\section{The classification theorem}\label{S:class}

   \begin{theorem} \label{teo-alg}
Let $\n$ be a complex $4$-dimensional nilpotent  algebra.
Then $\n$ is isomorphic to an algebra from the following list:
 
{\tiny
\begin{longtable}{lllll llll}
\hline$\cd {4*}{01}$&$:$& $e_1 e_1 = e_2$\\
\hline$\cd {4*}{02}$&$:$& $e_1 e_1 = e_3$ &$ e_2 e_2=e_3$ \\
\hline$\cd {4*}{03}$&$:$& $e_1 e_2=e_3$ & $e_2 e_1=-e_3$ \\
\hline$\cd {4*}{04}(\lambda)$&$:$& $e_1 e_1 = \lambda e_3$ & $e_2 e_1=e_3$  & $e_2 e_2=e_3$\\

\hline$\cd {4*}{05}$ & $:$ & $e_1e_1 = e_3$&$ e_2e_2=e_4$ \\ 
\hline$\cd {4*}{06}$ & $:$ & $e_1e_2 = e_4$&$ e_3e_1 = e_4$   \\ 
\hline$\cd {4*}{07}$ & $:$ & $e_1e_2 = e_3$&$ e_2e_1 = e_4$&$  e_2e_2 = -e_3$\\ 
\hline$\cd {4*}{08}(\af)$ & $:$&$ e_1e_1 = e_3$& $e_1e_2 = e_4$&  $e_2e_1 = -\af e_3$&$e_2e_2 = -e_4$ \\
\hline$\cd {4*}{09}(\alpha)$&$:$&$e_1e_1 = e_4$&$ e_1e_2 = \alpha e_4$&$  e_2e_1 = -\alpha e_4$\\&&$ e_2e_2 = e_4$&$  e_3e_3 = e_4$\\
\hline$\cd {4*}{10}$&$:$&$ e_1e_2 = e_4$&$ e_1e_3 = e_4$&$ e_2e_1 = -e_4$\\&&$ e_2e_2 = e_4$&$ e_3e_1 = e_4$\\ 
\hline$\cd {4*}{11}$&$:$&$ e_1e_1 = e_4$&$ e_1e_2 = e_4$&$ e_2e_1 = -e_4$&$ e_3e_3 = e_4$\\
\hline$\cd {4*}{12}$&$:$&$ e_1e_2 = e_3$&$ e_2e_1 = e_4$  \\ 
\hline$\cd {4*}{13}$&$:$&$ e_1e_1 = e_4$&$ e_1e_2 = e_3$&$ e_2e_1 = -e_3$&$ e_2e_2=2e_3+e_4$\\
\hline$\cd {4*}{14}(\af)$&$:$&$ e_1e_2 = e_4$&$ e_2e_1 =\alpha e_4$&$ e_2e_2 = e_3$\\
\hline$\cd {4*}{15}$&$:$&$e_1e_2 = e_4$&$ e_2e_1 = -e_4$&$ e_3e_3 = e_4$\\
\hline

$\D{4}{00} $&$:$&  
$e_1e_1=e_4$ & $e_1e_2=e_4$& $e_2 e_1=e_3$ \\&& $e_2 e_2 = e_3$ & $e_2e_3=e_4$ 
\\\hline

$\D{4}{01}(\lambda,\alpha,\beta)$&$:$& 
$e_1 e_1 = \lambda e_3 + e_4$ & $e_1 e_3 = \alpha e_4$ & $e_2 e_1=e_3$ \\
&&$ e_2 e_2 = e_3$ & $e_2 e_3 = \beta e_4$ & $e_3e_1 = e_4$ \\
\hline
$\D{4}{02}(\lambda,\alpha,\beta)$&$:$& 
$e_1 e_1 = \lambda e_3$ & $e_1 e_3 = \alpha e_4$ & $e_2 e_1=e_3$ \\
&& $e_2 e_2 = e_3$ & $e_2 e_3 = \beta e_4$ & $e_3e_1 = e_4$ \\
\hline
$\D{4}{03}(\lambda,\alpha)$&$:$& 
$e_1 e_1 = \lambda e_3$ & $e_1 e_2 = e_4$ & $e_2 e_1=e_3$ \\
&& $e_2 e_2 = e_3$ & $e_2 e_3 = \alpha e_4$ & $e_3e_1 = e_4$ \\
\hline
$\D{4}{04}(\lambda,\alpha)$&$:$& 
$e_1 e_1 = \lambda e_3$ & $e_2 e_1=e_3$ & $e_2 e_2 = e_3$ \\
&& $e_2 e_3 = \alpha e_4$ & $e_3e_1 = e_4$ \\
\hline
$\D{4}{05}(\lambda,\alpha)$&$:$& 
$e_1 e_1 = \lambda e_3$ & $e_1 e_2 = \lambda e_4$ & $e_2 e_1=e_3 + \lambda\alpha e_4$ \\
&& $e_2 e_2 = e_3$ & $e_2 e_3 = 	\Theta e_4$ & $e_3e_1 = \lambda e_4$ \\
\hline
$\D{4}{06}(\lambda,\alpha)$&$:$& 
$e_1 e_1 = \lambda e_3$ & $e_1 e_2 = e_4$ & $e_2 e_1=e_3 + \alpha e_4$ \\
&& $e_2 e_2 = e_3$ & $e_2 e_3 = 	\Theta^{-1} e_4$ & $e_3e_1 = e_4$ \\
\hline
$\D{4}{07}(\lambda)$&$:$& $e_1 e_1 = \lambda e_3$ & $e_2 e_1=e_3 + \lambda e_4$ & $e_2 e_2 = e_3$ \\&& $e_2 e_3 = 	\Theta e_4$ & $e_3e_1 = \lambda e_4$ \\
\hline
$\D{4}{08}(\lambda)$&$:$& 
$e_1 e_1 = \lambda e_3$ &  $e_2 e_1=e_3 + e_4$ & $e_2 e_2 = e_3$ \\&& $e_2 e_3 = 	\Theta^{-1} e_4$ & $e_3e_1 = e_4$ \\
\hline
$\D{4}{09}(\lambda,\alpha)$&$:$& 
$e_1 e_1 = \lambda e_3$ & $e_1 e_2 = e_4$& $e_1 e_3 = \alpha e_4$ \\
&& $e_2 e_1=e_3$ & $e_2 e_2 = e_3 $& $e_3e_1 = e_4$ \\
\hline
$\D{4}{10}(\lambda,\alpha)$&$:$& 
$e_1 e_1 = \lambda e_3$& $e_1 e_3 = \alpha e_4$ &  $e_2 e_1=e_3$ \\&& $e_2 e_2 = e_3$ & $e_3e_1 = e_4$ \\
\hline        
$\D{4}{11}(\lambda,\alpha)$&$:$& 
$e_1 e_1 = \lambda e_3 + e_4$& $e_1e_2 = \alpha e_4$ & $e_1 e_3 = -e_4$ \\
&&  $e_2 e_1=e_3$ & $e_2 e_2 = e_3$ & $e_3e_1 = e_4$ \\
\hline
$\D{4}{12}(\lambda,\alpha)$&$:$& 
$e_1 e_1 = \lambda e_3 + e_4$& $e_1e_3 = \alpha e_4$  &  $e_2 e_1=e_3$ \\
&& $e_2 e_2 = e_3$ & $e_3e_1 = \Theta e_4$ & $e_3e_2 = e_4$ \\
\hline
$\D{4}{13}(\lambda,\alpha)$&$:$& 
$e_1 e_1 = \lambda e_3 + e_4$& $e_1e_3 = \alpha e_4$  &  $e_2 e_1=e_3$ \\
&& $e_2 e_2 = e_3$ & $e_3e_1 = (1-\Theta)e_4$ & $e_3e_2 = e_4$ \\
\hline
$\D{4}{14}(\lambda,\alpha)$&$:$& 
$e_1 e_1 = \lambda e_3$& $e_1e_3 = \alpha e_4 $ &  $e_2 e_1=e_3$ \\
&& $e_2 e_2 = e_3$ & $e_3e_1 = \Theta e_4$ & $e_3e_2 = e_4$ \\
\hline
$\D{4}{15}(\lambda,\alpha)$&$:$& 
$e_1 e_1 = \lambda e_3$& $e_1e_3 = \alpha e_4$  &  $e_2 e_1=e_3$ \\
&& $e_2 e_2 = e_3$ & $e_3e_1 = (1-\Theta)e_4$ & $e_3e_2 = e_4$ \\
\hline
$\D{4}{16}(\alpha)$&$:$& $e_1e_3 = \alpha e_4$  &  $e_2 e_1=e_3 + e_4$ & $e_2 e_2 = e_3$ \\&& $e_3e_1 = e_4$ & $e_3e_2 = e_4$ \\
\hline
$\D{4}{17}(\alpha)$&$:$& 
$e_1 e_1 = e_4$ & $e_1e_3 = -e_4$  &  $e_2 e_1=e_3 + \alpha e_4$ \\
&& $e_2 e_2 = e_3$ & $e_3e_1 = e_4$ & $e_3e_2 = e_4$ \\
\hline
$\D{4}{18}(\lambda,\alpha)$&$:$& $e_1 e_1 = \lambda e_3 + e_4$& $e_1e_3 = \Theta\alpha e_4$  &  $e_2 e_1=e_3$ & $e_2 e_2 = e_3$ \\ && $e_2 e_3 = \alpha e_4$ & $e_3e_1 = \Theta e_4$ & $e_3e_2 = e_4$ \\

\hline
$\D{4}{19}(\lambda,\alpha)$&$:$& $e_1 e_1 = \lambda e_3 + e_4$& $e_1e_3 = (1-\Theta)\alpha e_4$  &  $e_2 e_1=e_3$ & $e_2 e_2 = e_3$ \\&& $e_2 e_3 = \alpha e_4$ & $e_3e_1 = (1-\Theta)e_4$& $e_3e_2 = e_4$ \\
\hline   
$\D{4}{20}(\lambda,\alpha)$&$:$& $e_1 e_1 = \lambda e_3$& $e_1e_3 = \Theta\alpha e_4$  &  $e_2 e_1=e_3$ & $e_2 e_2 = e_3$ \\
&& $e_2 e_3 = \alpha e_4$ & $e_3e_1 = \Theta e_4$& $e_3e_2 = e_4$ \\
\hline
$\D{4}{21}(\lambda,\alpha)$&$:$& $e_1 e_1 = \lambda e_3$& $e_1e_3 = (1-\Theta)\alpha e_4$  &  $e_2 e_1=e_3$ & $e_2 e_2 = e_3$ \\&& $e_2 e_3 = \alpha e_4$ & $e_3e_1 = (1-\Theta)e_4$ & $e_3e_2 = e_4$ \\
\hline
$\D{4}{22}(\lambda)$&$:$& $e_1 e_1 = \lambda e_3 + (1-2\lambda)e_4$& $e_1 e_2 = e_4$& $e_1e_3 = (\Theta - 1) e_4$  &  $e_2 e_1=e_3$ \\&& $e_2 e_2 = e_3$ & $e_2 e_3 = (1-\Theta^{-1}) e_4$& $e_3e_1 = \Theta e_4$ & $e_3e_2 = e_4$ \\
\hline
$\D{4}{23}(\lambda)$&$:$& $e_1 e_1 = \lambda e_3 + \lambda(1-2\lambda) e_4$& $e_1 e_2 = \lambda e_4$& $e_1e_3 = -\lambda\Theta e_4 $&  $e_2 e_1=e_3$ \\&& 
$e_2 e_2 = e_3$ & $e_2 e_3 = -\Theta^2 e_4$ & $e_3e_1 = \lambda(1-\Theta)e_4$ & $e_3e_2 = \lambda e_4$ \\
\hline        
$\D{4}{24}(\lambda)$&$:$& $e_1 e_1 = \lambda e_3 + \Theta e_4$& $e_1 e_2 = e_4$& $e_1e_3 = (\Theta - 1) e_4$  &  $e_2 e_1=e_3$ \\&& $e_2 e_2 = e_3$ & $e_2 e_3 = (1-\Theta^{-1}) e_4$ & $e_3e_1 = \Theta e_4$ & $e_3e_2 = e_4$ \\
\hline
$\D{4}{25}(\lambda)$&$:$& $e_1 e_1 = \lambda e_3 + \lambda(1-\Theta)e_4$& $e_1 e_2 = \lambda e_4$& $e_1e_3 = -\lambda\Theta e_4$  &  $e_2 e_1=e_3$ \\&& $e_2 e_2 = e_3$ & $e_2 e_3 = -\Theta^2 e_4$ 
        & $e_3e_1 = \lambda (1-\Theta)e_4$ & $e_3e_2 = \lambda e_4$ \\
\hline
$\D{4}{26}(\lambda)$&$:$& $e_1 e_1 = \lambda e_3 + \Theta e_4$& $e_1 e_2 = e_4$& $e_1e_3 = -\Theta e_4$  &  $e_2 e_1=e_3$\\
& & $e_2 e_2 = e_3$ & $e_2 e_3 = -e_4$
         & $e_3e_1 = \Theta e_4$ & $e_3e_2 = e_4 $\\
\hline
$\D{4}{27}(\lambda)$&$:$& $e_1 e_1 = \lambda e_3 + (1-\Theta)e_4$& $e_1 e_2 = e_4$& $e_1e_3 = (\Theta-1) e_4$  &  $e_2 e_1=e_3$ \\&& $e_2 e_2 = e_3$ & $e_2 e_3 = -e_4$  & $e_3e_1 = (1-\Theta)e_4$ & $e_3e_2 = e_4$ \\
\hline
$\D{4}{28}(\lambda)$&$:$& $e_1 e_1 = \lambda e_3 + (1-\Theta)e_4$& $e_1 e_2 = e_4$& $e_1e_3 = -\Theta e_4$  &  $e_2 e_1=e_3$ \\&& $e_2 e_2 = e_3$ & $e_2 e_3 = -e_4$ & $e_3e_1 = \Theta e_4$ & $e_3e_2 = e_4$ \\
\hline
$\D{4}{29}(\lambda)$&$:$& $e_1 e_1 = \lambda e_3 + \Theta e_4$& $e_1 e_2 = e_4$& $e_1e_3 = (\Theta-1) e_4$  &  $e_2 e_1=e_3$ \\
&& $e_2 e_2 = e_3$ & $e_2 e_3 = -e_4$ & $e_3e_1 = (1-\Theta)e_4$ & $e_3e_2 = e_4$ \\
$\D{4}{30}(\lambda)$&$:$& $e_1 e_1 = \lambda e_3 + e_4$ & $e_1e_3 = (\Theta-1) e_4$  &  $e_2 e_1=e_3$ & $e_2 e_2 = e_3$ \\
& & $e_2 e_3 = -e_4$ & $e_3e_1 = \Theta e_4$ & $e_3e_2 = e_4$ \\
\hline
$\D{4}{31}(\lambda)$&$:$& $e_1 e_1 = \lambda e_3 + e_4$ & $e_1e_3 = -\Theta e_4$  &  $e_2 e_1=e_3$ & $e_2 e_2 = e_3$ \\
&& $e_2 e_3 = -e_4$ & $e_3e_1 = (1-\Theta)e_4$ & $e_3e_2 = e_4$ \\
\hline        
$\D{4}{32}(\lambda)$&$:$& $e_1 e_1 = \lambda e_3$ & $e_1e_3 = (\Theta-1) e_4$  &  $e_2 e_1=e_3$ & $e_2 e_2 = e_3$ \\
&& $e_2 e_3 = -e_4$ & $e_3e_1 = \Theta e_4$ & $e_3e_2 = e_4$ \\
\hline
$\D{4}{33}(\lambda)$&$:$& $e_1 e_1 = \lambda e_3$ & $e_1e_3 = -\Theta e_4$  & $ e_2 e_1=e_3$ & $e_2 e_2 = e_3$ \\
&& $e_2 e_3 = -e_4$ & $e_3e_1 = (1-\Theta)e_4$ & $e_3e_2 = e_4$ \\
\hline 
$\D{4}{34}$&$:$& $e_1 e_1 = e_4$ &  $e_2 e_1=e_3 + e_4$ & $e_2 e_2 = e_3$\\& & $e_3 e_1 = e_4$& $e_3 e_2 = e_4$ \\
\hline
$\D{4}{35}(\lambda)$&$:$& $e_1 e_1 = \lambda e_3$ & $e_1e_3 = e_4$  &  $e_2 e_1=e_3 + e_4$ & $e_2 e_2 = e_3$ \\
\hline
$\D{4}{36}(\lambda)$&$:$& $e_1 e_1 = \lambda e_3$ & $e_1e_3 = e_4$  &  $e_2 e_1=e_3$ & $e_2 e_2 = e_3$ \\
\hline 
$\D{4}{37}(\lambda)$&$:$& $e_1 e_1 = \lambda e_3 + e_4$ & $e_1e_3 = \Theta e_4$ &  $e_2 e_1=e_3$ \\&& $e_2 e_2 = e_3$&  $e_2 e_3 = e_4$ \\
\hline
$\D{4}{38}(\lambda)$&$:$& $e_1 e_1 = \lambda e_3 + e_4$ & $e_1e_3 = (1-\Theta)e_4$ &  $e_2 e_1=e_3$ \\&& $e_2 e_2 = e_3$&  $e_2 e_3 = e_4$ \\
\hline
$\D{4}{39}(\lambda)$&$:$& $e_1 e_1 = \lambda e_3$ & $e_1e_3 = \Theta e_4$ &  $e_2 e_1=e_3$ \\&& $e_2 e_2 = e_3$&  $e_2 e_3 = e_4$ \\
\hline
$\D{4}{40}(\lambda)$&$:$& $e_1 e_1 = \lambda e_3$ & $e_1e_3 = (1-\Theta)e_4$ &  $e_2 e_1=e_3$ \\&& $e_2 e_2 = e_3$&  $e_2 e_3 = e_4$\\
\hline
$\cd 4{01}$&$:$& $e_1 e_1 = e_2$ & $e_2 e_2=e_3$ \\
\hline
$\cd 4{02}$&$:$& $e_1 e_1 = e_2$ & $e_2 e_1= e_3$ & $e_2 e_2=e_3$ \\
\hline$\cd 4{03}$&$:$& $e_1 e_1 = e_2$ & $e_2 e_1=e_3$ \\
\hline$\cd 4{04}(\lambda)$&$:$& $ e_1 e_1 = e_2$ & $e_1 e_2=e_3$ & $e_2 e_1=\lambda e_3$ \\ \hline

$\cd 4{05}$ & $:$ & $e_1 e_1 = e_2$ & $e_2 e_1=e_4$ & $e_2 e_2=e_3$ \\
\hline$\cd 4{06}$ & $:$ & $ e_1 e_1 = e_2$ & $e_1 e_2=e_4$ & $e_2 e_1=e_3$  \\

\hline$\cd 4{07}(\lambda)$& $:$ & $e_1 e_1 = e_2$ & $e_1 e_2=e_4$ & $e_2 e_1=\lambda e_4$ & $e_2 e_2=e_3$ \\\hline
$\cd 4{08}(\alpha)$ &$:$&
$e_1 e_1 = e_2$ & $e_1e_3=e_4$ & $e_2 e_1=e_3$ \\&& $e_2e_2=\alpha e_4$ & $e_3e_1=-2e_4$ \\

\hline$\cd 4{09}$ &$:$& 
$e_1 e_1 = e_2$ & $e_1e_2=e_4$ & $e_1e_3=e_4$ \\
&& $e_2 e_1=e_3$ & $e_2e_2=- e_4$ & $e_3e_1=-2e_4$\\
\hline

$\cd 4{10}(\alpha)$ &$:$&
$e_1 e_1 = e_2$ & $e_1 e_2=e_3$ & $e_1e_3=-e_4$ \\
&& $e_2 e_1=  e_3 +e_4$ & $e_2e_2=\alpha e_4$ & $e_3e_1=- e_4$ \\

\hline$\cd 4{11}(\lambda)$ &$:$& 
$e_1 e_1 = e_2$ & $e_1 e_2=e_3$ & $e_1e_3=(\lambda-2)e_4$ \\
$\lambda \neq1$&& $e_2 e_1=  \lb e_3+e_4$  &
$e_2e_2=-(\lambda+1)^2 e_4$ &$e_3e_1= (1-2\lb) e_4$ \\

\hline$\cd 4{12}(\alpha, \lambda)$ &$:$& 
$e_1 e_1 = e_2$ & $e_1 e_2=e_3$ & $e_1e_3=(\lambda-2)e_4$ \\&& $e_2 e_1=\lambda e_3$ & $e_2e_2=\alpha e_4$ & $e_3e_1=(1-2\lambda) e_4$ \\
\hline

$\cd {4}{13}(\af)$&$:$& 
$e_1 e_1 = e_2$ & $e_1e_2=e_4$ & $e_1e_3=e_4$& $e_2e_1=e_4$ \\
$\af\neq \frac{1}{2}$&&  $e_2e_3=\alpha e_4$& $e_3e_1=e_4$& \multicolumn{2}{l}{$e_3e_2=(\alpha +1)e_4$}\\

\hline
$\cd {4}{14}(\af, \beta)$&$:$& 
$e_1 e_1 = e_2$ & $e_1 e_2 = e_4$&  $e_1 e_3 = \af e_4$& $e_2 e_1 = e_4$\\ && 
$e_2 e_2 = e_4$& $e_3 e_1 = \af e_4$& $e_3 e_2 = e_4$& $e_3 e_3 =\beta  e_4$& \\

\hline
$\cd {4}{15}(\af)$&$:$& 
$e_1 e_1 = e_2$ &  $e_1 e_2 = \af e_4$&   $e_1 e_3 = e_4$\\&&
 $e_2 e_1 =(\af+1)  e_4$ &  $e_3 e_1 = e_4$\\

\hline
$\cd {4}{16}$&$:$& 
$e_1 e_1 = e_2$ & $e_1e_2=e_4$ & $e_2 e_1 = e_4$\\ && $e_2 e_3 = -\frac{1}{2} e_4$& $e_3 e_2 =\frac{1}{2} e_4$& $e_3 e_3 = e_4$&\\

\hline
$\cd {4}{17}(\af)$&$:$& 
$e_1 e_1 = e_2$ &  $e_1 e_2 = e_4$&  $e_2 e_1 = e_4$\\&&  $e_2 e_3 =\af e_4$&  $e_3 e_2 = (\af+1) e_4$&\\

\hline
$\cd {4}{18}(\af)$&$:$& 
$e_1 e_1 = e_2$ &  $e_1 e_2 =\af  e_4$&   \multicolumn{1}{l}{$e_2 e_1 =(\af+1) e_4$}&   $e_3 e_3 = e_4$&\\

\hline
$\cd {4}{19}$&$:$& 
$e_1 e_1 = e_2$ &   $e_1 e_2 = e_4$&   $e_2 e_1 = e_4$\\&&    $e_3 e_1 = e_4$&
 $e_3 e_3 = e_4$&\\
 
\hline
$\cd {4}{20}$&$:$& 
$e_1 e_1 = e_2$ &    $e_1 e_2 = e_4$&   $e_2 e_1 = e_4$&   $e_3 e_1 = e_4$&\\

\hline
$\cd {4}{21}(\af)$&$:$& 
$e_1 e_1 = e_2$ & $e_1 e_3 =\af e_4$&   $e_2 e_1 = e_4$&  $e_2 e_2 = e_4$&\\
&&  $e_2 e_3 = e_4$&  $e_3 e_1 =\af e_4$&  $e_3 e_2 = e_4$&  $e_3 e_3 = e_4$\\

\hline
$\cd {4}{22}$&$:$& 
$e_1 e_1 = e_2$ &   $e_1 e_3 = e_4$&    $e_2 e_3 = -\frac{1}{2} e_4$\\&&
 $e_3 e_1= e_4$&   $e_3 e_2 =\frac{1}{2} e_4$&   $e_3 e_3 = e_4$&\\

\hline
$\cd {4}{23}(\af)$&$:$& 
  $e_1 e_1 = e_2$&   $e_1 e_3 = e_4$&   $e_2 e_3 = \af e_4$\\&&
  $e_3 e_1 = e_4$&   $e_3 e_2 =(\af +1)e_4$&\\

\hline
$\cd {4}{24}(\af)$&$:$& 
  $e_1 e_1 = e_2$&   $e_1 e_3 = e_4$&    $e_2 e_2 = e_4$\\
  && $e_3 e_1 = e_4$&   $e_3 e_2 = e_4$& $e_3 e_3 =\af  e_4$&\\

\hline
$\cd {4}{25} $&$:$& 
  $e_1 e_1 = e_2$&  $e_1 e_3 = e_4$& $e_2 e_1 = e_4$\\&& $e_3 e_1 = e_4$&    $e_2 e_2 = e_4$&\\

\hline
$\cd {4}{26}(\af)$&$:$& 
  $e_1 e_1 = e_2$&   $e_1 e_3 = \af e_4$&    $e_2 e_2 = e_4$& \multicolumn{2}{l}{$e_3 e_1 =(\af+1) e_4$}&\\

  
  \hline
$\cd {4}{27}$&$:$& 
  $e_1 e_1 = e_2$&   $e_2 e_1 = e_4$& $e_2 e_3 = e_4$\\&&  $e_3 e_2 = e_4$&   $e_3 e_3 = e_4$&\\

  \hline
$\cd {4}{28}(\af)$&$:$& 
  $e_1 e_1 = e_2$&   $e_2 e_1 = e_4$&    $e_2 e_2 = e_4$\\
  $\af\ne 1$&& $e_2 e_3 =   e_4$&   $e_3 e_2 =  e_4$&   $e_3 e_3 =\af  e_4$&\\

  \hline
$\cd {4}{29}$&$:$& 
  $e_1 e_1 = e_2$&    $e_2 e_3 =-\frac{1}{2} e_4$&   $e_3 e_2 =\frac{1}{2} e_4$&   $e_3 e_3 = e_4$&\\

  \hline
$\cd {4}{30}$&$:$& 
  $e_1 e_1 = e_2$&   $e_2 e_1 = e_4$& $e_2 e_3 = e_4$&   $e_3 e_2 = e_4$&\\

  \hline
$\cd {4}{31}$&$:$& 
  $e_1 e_1 = e_2$&  $e_2 e_3 = e_4$&     $e_3 e_1 = e_4$&
$e_3 e_2 = e_4$&\\

  \hline
$\cd {4}{32}(\af)$&$:$& 
  $e_1 e_1 = e_2$&   $e_2 e_3 = \af e_4$&   \multicolumn{2}{l}{$e_3 e_2 = (\af+1) e_4$}& \\

  \hline
$\cd {4}{33}$&$:$& 
  $e_1 e_1 = e_2$&   $e_2 e_1 = e_4$&   $e_2 e_2 = e_4$& $e_3 e_3 = e_4$& \\

 \hline
$\cd {4}{34}$&$:$& 
  $e_1 e_1 = e_2$&    $e_2 e_2 = e_4$& $e_3 e_3 = e_4$& \\

 \hline
$\cd {4}{35}$&$:$& 
  $e_1 e_1 = e_2$&    $e_2 e_2 = e_4$&  $e_3 e_1 = e_4$& $e_3 e_3 = e_4$& \\

 \hline
$\cd {4}{36}(\af)$&$:$& 
  $e_1 e_1 = e_2$&     $e_2 e_2 = e_4$&   $e_3 e_2= e_4$& $e_3 e_3 =\af  e_4$&\\


 \hline
$\cd {4}{37}$&$:$& 
  $e_1 e_1= e_2$&    $e_1 e_2= e_4$&   $e_2 e_1= e_4$&    $e_3 e_3= e_4$\\


  \hline
$\cd {4}{38}$&$:$& 
  $e_1 e_1= e_2$&   $e_2 e_3= e_4$&   $e_3 e_2= e_4$&\\

  
\hline
$\cd 4{39}$ & $:$ &  
$e_1 e_1 = e_3+e_4$ & $e_1e_2=\frac i2 e_4$ & $e_1e_3=e_4$ & $e_2e_1=\frac i2 e_4$\\
&& $ e_2 e_2=e_3$ & $e_2e_3=-2ie_4$ & $e_3e_1=2e_4$ & $e_3e_2=-ie_4$\\

\hline
$\cd 4{40}$ & $:$ &  
$e_1 e_1 = e_3+e_4$ & $e_1e_2=\frac i2 e_4$ & $e_1e_3=-\frac 12e_4$ & $e_2e_1=\frac i2 e_4$\\
&& $ e_2 e_2=e_3$ & $e_2e_3=-\frac i2e_4$ & $e_3e_1=\frac 12e_4$ & $e_3e_2=\frac i2e_4$\\

\hline
$\cd 4{41}$ & $:$ &  
$e_1 e_1 = e_3+e_4$ & $e_1e_2=-\frac i2 e_4$ & $e_1e_3=e_4$ & $e_2e_1=-\frac i2 e_4$\\
&& $ e_2 e_2=e_3$ & $e_2e_3=-2ie_4$ & $e_3e_1=2e_4$ & $e_3e_2=-ie_4$\\

\hline
 $\cd 4{42}$ & $:$ &  
$e_1 e_1 = e_3+e_4$ & $e_1e_2=-\frac i2 e_4$ & $e_1e_3=-\frac 12e_4$ & $e_2e_1=-\frac i2 e_4$\\
&& $ e_2 e_2=e_3$ & $e_2e_3=-\frac i2e_4$ & $e_3e_1=\frac 12e_4$ & $e_3e_2=\frac i2e_4$\\

\hline
$\cd 4{43}(\af)$ &$:$&  
$e_1 e_1 = e_3 + e_4$     & $e_1 e_2 = \af e_4$          & $e_1 e_3 = -\frac 12 e_4$ \\
&&   $e_2 e_1 = \af e_4$      & $e_2 e_2 = e_3$           & $e_3 e_1 = \frac 12 e_4$\\

\hline
$\cd 4{44}(\af,\bt,\gm)$ & $:$ &  
$e_1 e_1 = e_3 + \af e_4$   & $e_1 e_2 = \bt e_4$ & $e_2 e_1 = (\bt+\gm) e_4$ \\
&& $e_2 e_2 = e_3$ & $e_3 e_1 = e_4$             & $e_3 e_3 = e_4$\\

\hline
$\cd 4{45}$ & $:$ &  
$e_1 e_1 = e_3 + 2i e_4$  & $e_1 e_2 = e_4$ & $e_2 e_1 = e_4$& $e_2 e_2 = e_3$\\
&&  $e_3 e_1 = e_4$           & $e_3 e_2 = i e_4$       & $e_3 e_3 = e_4$\\

\hline
$\cd 4{46}(\af)$ & $:$ &
$e_1 e_1 = e_3 - 2i\af e_4$ & $e_1 e_2 = \af e_4$  & $e_2 e_1 = \af e_4$ & $e_2 e_2 = e_3$ \\
$\af\ne 0$& & $e_3 e_1 = e_4$              & $e_3 e_2 = i e_4$       & $e_3 e_3 = e_4$\\

\hline
$\cd 4{47}(\af,\bt)$ & $:$ & 
$e_1 e_1 = e_3 + e_4$  & $e_1 e_3 = \af e_4$  & $e_2 e_2 = e_3$\\ 
$\bt\ne 0$ && $e_2 e_3 = \bt e_4$  & $e_3 e_1 =(\af+1) e_4$        & $e_3 e_2 = \bt e_4$ \\

\hline		
$\cd 4{48}(\af)$ & $:$ & 
 $e_1 e_1 = e_3 + \af e_4$ &  $e_2 e_1 = i\af e_4$ & $e_2 e_2 = e_3$  \\
$\af\ne 0$&& $e_3 e_1 = e_4$            & $e_3 e_2 = i e_4$     & $e_3 e_3 = e_4$\\

\hline 
$\cd 4{49}(\af)$ & $:$ & 
$e_1 e_1 = e_3 + \af e_4$ & $e_2 e_1 = -i\af e_4$ & $e_2 e_2 = e_3$\\
$\af\ne 0$& & $e_3 e_1 = e_4$& $e_3 e_2 = i e_4$ & $e_3 e_3 = e_4$\\

\hline 
$\cd 4{50}(\af)$ & $:$ & 
$e_1 e_1 = e_3 + \af e_4$ & $e_2 e_1 = e_4$ & $e_2 e_2 = e_3$ & $e_3 e_3 = e_4$\\

\hline
$\cd 4{51}(\af)$ & $:$ &
$e_1 e_1 = e_3 + \af e_4$ &   $e_2 e_2 = e_3$ & $e_3 e_1 = e_4$\\
&& $e_3 e_2 = i e_4$ & $e_3 e_3 = e_4$\\

\hline 
$\cd 4{52}$ & $:$ &  
$e_1 e_1 = e_3 +  e_4$ &   $e_2 e_2 = e_3$ &   $e_3 e_3 = e_4$\\

\hline 
$  \cd 4{53}$ & $:$ &  
$e_1 e_1 = e_3$ &  $e_1e_2=-\frac 12e_4$ & $e_1e_3=e_4$ & $e_2e_1=\frac 12e_4$\\
 && $e_2 e_2 = e_3$ & $e_2e_3=ie_4$ &  $e_3e_1=e_4$ & $e_3e_2=ie_4$ \\

\hline
 $\cd 4{54}(\af)$ & $:$ &  
$e_1 e_1 = e_3$ & $e_1e_2=e_4$ & $e_1e_3=\af e_4$ & $e_2e_1=e_4$\\
&& $ e_2 e_2=e_3$ & $e_2e_3=-i(\af+1)e_4$ & $e_3e_1=(\af+1) e_4$ & $e_3e_2=-i\af e_4$\\

\hline
$\cd 4{55}(\af)$ & $:$ &  
$e_1 e_1 = e_3$ & $e_1 e_2 = e_4$ & $e_1 e_3 = \af e_4$ \\
& & $e_2 e_1 = e_4$ & $e_2 e_2 = e_3$  & $e_3 e_1 = (\af+1) e_4$ \\
		
\hline
$\cd 4{56}$ & $:$ &  
$e_1 e_1 = e_3$ & $e_1 e_2 = e_4$ & $e_2 e_1 = -e_4$ & $e_2 e_2 = e_3$ \\
& & $e_3 e_1 = e_4$ & $e_3 e_2 = i e_4$ & $e_3 e_3 = e_4$\\

\hline
$\cd 4{57}(\af,\bt)$ & $:$ &  
$e_1 e_1 = e_3$ & $e_1 e_2 = \af e_4$ &   $e_2 e_1 = (\af+\bt)e_4$ & $e_2 e_2 = e_3$ \\
$\bt\not\in\{0,-2\af\}$& & $e_3 e_1 = e_4$          & $e_3 e_2 = i e_4$     & $e_3 e_3 = e_4$\\

\hline
$\cd 4{58}$ & $:$ &  
$e_1 e_1 = e_3$ & $e_1 e_3 = e_4$  & $e_2 e_1 = e_4$ & $e_2 e_2 = e_3$\\
&& $e_2 e_3 = i e_4$  & $e_3 e_1 = e_4$ & $e_3 e_2 = i e_4$ \\

\hline 
$\cd 4{59}(\af,\bt)$ & $:$ &   
$e_1 e_1 = e_3$ & $e_1 e_3 = \af e_4$ & $e_2 e_2 = e_3$\\
 $\bt\ne 0$ && $e_2 e_3 = \bt e_4$ & $e_3 e_1 = (\af+1) e_4$  & $e_3 e_2 = \bt e_4$ \\

\hline		
$\cd 4{60}$ & $:$ & 
$e_1 e_1 = e_3$ & $e_1 e_3 = i e_4$  & $e_2 e_2 = e_3$\\
&& $e_2 e_3 = e_4$ & $e_3 e_1 = (i+1) e_4$ & $e_3 e_2 = (i+1) e_4$  \\

\hline
$\cd 4{61}(\af)$ & $:$ & 
 $e_1 e_1 = e_3$ & $e_1 e_3 = -i\af e_4$ & $e_2 e_2 = e_3$\\
 && $e_2 e_3 = \af e_4$ & $e_3 e_1 = (1-i\af) e_4$ & $e_3 e_2 = (\af+i) e_4$  \\

\hline
$\cd 4{62}$ & $:$ &  
$e_1 e_1 = e_3$  & $e_1e_3=e_4$ & $e_2 e_2=e_3$\\
&& $e_2e_3=-2ie_4$ & $e_3e_1=2e_4$ & $e_3e_2=-ie_4$\\

\hline
$\cd 4{63}$ & $:$ &  
$e_1 e_1 = e_3$ & $e_1e_3=-\frac 12e_4$ & $ e_2 e_2=e_3$\\
&& $e_2e_3=-\frac i2e_4$ & $e_3e_1=\frac 12e_4$ & $e_3e_2=\frac i2e_4$\\

\hline 
$\cd 4{64}(\af)$ & $:$ & 
$e_1 e_1 = e_3$ & $e_1 e_3 = \af e_4$ & $e_2 e_2 = e_3$  & $e_3 e_1 = (\af+1) e_4$ 		\\

\hline
$\cd 4{65}(\af)$ & $:$ &  
$e_1 e_1 = e_3$ &   $e_1 e_3 = \af e_4$ & $e_2 e_2 = e_3$ \\
$\af\ne 0$&&  $e_3 e_1 = (\af+1) e_4$  & $e_3 e_2 = i e_4$ \\
\hline
$\cd 4{66}$ & $:$ &  
$e_1 e_1 = e_3$& $e_2 e_1 = e_4$ & $e_2 e_2 = e_3$ \\
&& $e_2 e_3 = e_4$ & $e_3 e_2 = e_4$ \\

\hline 
$\cd 4{67}$ & $:$ & 
$e_1 e_1 = e_3$ &  $e_2 e_2 = e_3$ &  $e_3 e_3 = e_4$\\

\hline
$\cd 4{68}$ & $:$ & 
$e_1 e_1 = e_3+e_4$ & $e_1 e_3 = i e_4$  &$e_2 e_2 = e_3$ &\\
&&$e_2 e_3 = e_4$&$e_3 e_1 =   ie_4$ & $e_3 e_2 = e_4$  \\

\hline
$\cd 4{69}$ & $:$ & 
$e_1 e_1 = e_3$ & $e_1 e_3 =  ie_4$  &$e_2 e_2 = e_3$ &\\
&&$e_2 e_3 = e_4$&$e_3 e_1 =   ie_4$ & $e_3 e_2 = e_4$  \\

\hline
$\cd 4{70}$ & $:$ & 
$e_1 e_1 = e_3$ & $e_1 e_3 =    e_4$  &$e_2 e_2 = e_3$ & $e_3 e_1 = e_4$ &   \\
\hline
		$\cd 4{71}$ & $:$ & $e_1 e_1 = e_4$ & $e_1 e_2 = e_3$  & $e_1 e_3 = e_4$ \\&& $e_2 e_1 = -e_3$  & $e_3 e_1 = -e_4$\\
		\hline
		$\cd 4{72}$ & $:$ & $e_1 e_1 = e_4$ & $e_1 e_2 = e_3$ & $e_1 e_3 = e_4$\\
		&& $e_2 e_1 = -e_3$  & $e_2 e_2 = e_4$ & $e_3 e_1 = -e_4$\\
		\hline
		$\cd 4{73}$ & $:$ & $e_1 e_2 = e_3 + e_4$ & $e_1 e_3 = e_4$ & $e_2 e_1 = -e_3$  & $e_3 e_1 = -e_4$\\
		\hline
		$\cd 4{74}(\af)$ & $:$ & $e_1 e_2 = e_3$ & $e_1 e_3 = (\af+1)e_4$ & $e_2 e_1 = -e_3$ & $e_3 e_1 = -\af e_4$\\
		\hline
		$\cd 4{75}(\af)$ & $:$ & $e_1 e_2 = e_3$ & $e_1 e_3 = (\af+1)e_4$ & $e_2 e_1 = -e_3$ \\&& $e_2 e_2 = e_4$ & $e_3 e_1 = -\af e_4$\\        \hline
		$\cd 4{76}$ & $:$ & $e_1 e_2 = e_3$ & $e_1 e_3 = e_4$ & $e_2 e_1 = -e_3$ \\&& $e_2 e_2 = e_4$ & $e_3 e_1 = -e_4$\\        
		\hline
		$\cd 4{77}$ & $:$ & $e_1 e_2 = e_3$ & $e_1 e_3 = e_4$ & $e_2 e_1 = -e_3$ \\&& $e_2 e_3 = e_4$ & $e_3 e_2 = -e_4$\\
		\hline
		$\cd 4{78}$ & $:$ & $e_1 e_2 = e_3$ & $e_1 e_3 = e_4$ & $e_2 e_1 = -e_3$\\
		&& $e_2 e_2 = e_4$ & $e_2 e_3 = e_4$ & $e_3 e_2 = -e_4$\\
		\hline
		$\cd 4{79}(\af)$ & $:$ & $e_1 e_2 = e_3+\af e_4$ & $e_1 e_3 = e_4$ & $e_2 e_1 = -e_3$ \\&& $e_2 e_3 = e_4$ & $e_3 e_3 = e_4$\\
		\hline
		$\cd 4{80}$ & $:$ & $e_1 e_2 = e_3+e_4$ & $e_1 e_3 = e_4$ & $e_2 e_1 = -e_3$ & $e_3 e_3 = e_4$\\
		\hline
		$\cd 4{81}$ & $:$ & $e_1 e_2 = e_3+e_4$ & $e_2 e_1 = -e_3$ & $e_3 e_3 = e_4$\\
		\hline
		$\cd 4{82}$ & $:$ & $e_1 e_2 = e_3$ & $e_1 e_3 = e_4$ & $e_2 e_1 = -e_3$ \\&& $e_2 e_2 = e_4$ & $e_3 e_3 = e_4$\\
		\hline
		$\cd 4{83}( \af)$ & $:$ & $e_1 e_2 = e_3$ & $e_2 e_1 = -e_3$ & $e_2 e_2 = \af e_4$ \\
		$\af\ne 0$&& $e_2 e_3 = e_4$ & $e_3 e_3 = e_4$\\
		\hline
		$\cd 4{84}$ & $:$ & $e_1 e_2 = e_3$ & $e_2 e_1 = -e_3$ & $e_2 e_2 = e_4$ & $e_3 e_3 = e_4$\\
		\hline
		$\cd 4{85}$ & $:$ & $e_1 e_2 = e_3$ & $e_2 e_1 = -e_3$ & $e_3 e_3 = e_4$\\
		\hline
		$\cd 4{86}$ & $:$ & $e_1 e_2 = e_3$ & $e_2 e_1 = -e_3$ & $e_2 e_3 = e_4$  &  $e_3 e_2 = -e_4$\\
\hline		
		
$\cd {4}{87}(\lambda)$&$:$& 
$e_1 e_1 = \lambda e_3+(2 \Theta-1)e_4$& $e_1 e_2=e_4$ & $e_1e_3=e_4$& 
\multicolumn{2}{l}{$e_2 e_1=e_3-(1- \Theta)^2 \lb^{-1}e_4$}\\
$\lb \neq 0, \frac{1}{4}$&& $e_2 e_2=e_3$& $e_2e_3=\Theta\lb^{-1}e_4$ 
&$e_3e_3=e_4$\\

\hline$\cd {4}{88}(\lambda)$&$:$& 
$e_1 e_1 = \lambda e_3+(1-2 \Theta)e_4$& $e_1 e_2=e_4$ & $e_1e_3=e_4$& 
\multicolumn{2}{l}{$e_2 e_1=e_3- \Theta^2 \lb^{-1}e_4$}\\
$\lb \neq 0, \frac{1}{4}$&& $e_2 e_2=e_3$& $e_2e_3=\Theta\lb^{-1}e_4$ 
&$e_3e_3=e_4$\\

\hline$\cd {4}{89}(\lambda)$&$:$& 
$e_1 e_1 = \lambda e_3+(2 \Theta-1)e_4$& $e_1 e_2=e_4$ & $e_1e_3=e_4$& 
\multicolumn{2}{l}{$e_2 e_1=e_3-(1- \Theta)^2 \lb^{-1}e_4$}\\
$\lb \neq 0, \frac{1}{4}$&& $e_2 e_2=e_3$& $e_2e_3=(1-\Theta)\lb^{-1}e_4$ 
&$e_3e_3=e_4$\\

\hline$\cd {4}{90}(\lambda)$&$:$& 
$e_1 e_1 = \lambda e_3+(1-2 \Theta)e_4$& $e_1 e_2=e_4$ & $e_1e_3=e_4$& 
\multicolumn{2}{l}{$e_2 e_1=e_3- \Theta^2 \lb^{-1}e_4$}\\
$\lb \neq 0, \frac{1}{4}$&& $e_2 e_2=e_3$& $e_2e_3=(1-\Theta)\lb^{-1}e_4$ 
&$e_3e_3=e_4$\\

\hline$\cd {4}{91}(\lambda, \af)$&$:$& 
$e_1 e_1 = \lambda e_3+(2 \Theta-1)e_4$& $e_1 e_2=e_4$ & $e_1e_3=\af e_4$& \\
$\lb \neq 0, \frac{1}{4}$&& $e_2 e_1=e_3- (1-\Theta)^2 \lb^{-1}e_4$&
$e_2 e_2=e_3$&  $e_3e_3=e_4$\\

\hline$\cd {4}{92}(\lambda, \af)$&$:$& 
$e_1 e_1 = \lambda e_3+(1-2 \Theta)e_4$& $e_1 e_2=e_4$ & $e_1e_3=\af e_4$& \\
$\lb \neq 0, \frac{1}{4}$&& $e_2 e_1=e_3-  \Theta^2 \lb^{-1}e_4$&
$e_2 e_2=e_3$&  $e_3e_3=e_4$\\

\hline$\cd {4}{93}( \af)$&$:$& 
$e_1 e_1 = e_4$ & $e_1 e_2=e_4$ &$e_1e_3=\af e_4$ & $e_2e_1=e_3+e_4$\\
&&$ e_2 e_2=e_3$ &$e_2e_3=\af e_4$ & $e_3e_3=e_4$\\

\hline$\cd {4}{94}(\af, \bt)$&$:$& 
$e_1 e_1 = e_4$ & $e_1e_2=e_4$& $e_1e_3=\af e_4$&\\
$\af\neq0$&&$e_2 e_1=e_3+\beta e_4$   & $e_2 e_2=e_3$   &$e_3e_3=e_4$\\

\hline$\cd {4}{95}(\af)$&$:$& 
$e_1 e_1 =  e_4$ & $e_1e_2=e_4$ &$e_1e_3=\af e_4$ & $e_2 e_1=e_3$\\
&& $e_2 e_2=e_3$ & $e_2e_3=\af e_4$&$e_3e_3=e_4$\\

\hline$\cd {4}{96}(\af)$&$:$& 
$e_1 e_1 =  e_4$ & $e_1e_2=e_4$ & $e_2 e_1=e_3+\af e_4$\\
&& $e_2 e_2=e_3$ & $e_2e_3= e_4$&$e_3e_3=e_4$\\

\hline$\cd {4}{97}(\lambda)$&$:$& 
$e_1 e_1 = \lambda e_3$ & $e_1e_2=e_4$& $e_1e_3=\Theta e_4 $&$e_2 e_1=e_3-e_4$ \\
&& $e_2 e_2=e_3$ &$e_2e_3=e_4$ &$e_3e_3=e_4$\\

\hline$\cd {4}{98}(\lambda)$&$:$& 
$e_1 e_1 = \lambda e_3$ & $e_1e_2=e_4$& $e_1e_3=(1-\Theta) e_4 $&$e_2 e_1=e_3-e_4$ \\
$\lb \neq \frac{1}{4}$ && $e_2 e_2=e_3$ &$e_2e_3=e_4$ &$e_3e_3=e_4$\\

\hline$\cd {4}{99}(\af)$&$:$& 
$e_1e_2=e_4$& $e_1e_3=e_4$& $e_2 e_1=e_3-e_4$\\ 
$\af\neq1$ && $e_2 e_2=e_3$ &$e_2e_3=\af e_4$&$e_3e_3=e_4$\\

\hline$\cd {4}{100}(\af)$&$:$& 
$e_1 e_1 = \frac{1}{4} e_3$ & $e_1e_2=e_4$& $e_1e_3=\af e_4$& $e_2 e_1=e_3-e_4$ \\
$\af\notin\{0, \frac{1}{2}\}$&&  $e_2 e_2=e_3$ &$e_2e_3=2 \af e_4$ &$e_3e_3=e_4$\\

 \hline$\cd {4}{101}(\af, \bt)$&$:$& 
 $e_1e_2=e_4$& $e_1e_3=\af e_4$& $e_2 e_1=e_3$  \\
 && $e_2 e_2=e_3$&$e_2e_3=\bt e_4$ &$e_3e_3=e_4$\\
 
 \hline$\cd {4}{102}(\lb, \af)$&$:$& 
 $e_1 e_1 = \lambda e_3$ & $e_1e_2=e_4$& $e_2 e_1=e_3-e_4$ \\
 $\lb\neq0$&& $e_2 e_2=e_3$&$ e_2e_3=\af e_4$  &$e_3e_3=e_4$\\

\hline$\cd {4}{103}$&$:$& 
$e_1 e_2 = e_4$ & $e_2 e_1=e_3-e_4$  & $e_2 e_2=e_3$ &$e_3e_3=e_4$\\
 
 \hline$\cd {4}{104}$&$:$& 
 $e_1 e_3 =  e_4$ & $e_2 e_1=e_3+e_4$  & $e_2 e_2=e_3$ \\&&$e_2e_3=e_4$&$e_3e_3=e_4$\\
 
 \hline$\cd {4}{105}(\lambda, \af,\bt)$&$:$& 
 $e_1 e_1 = \lambda e_3$ & $e_1e_3=e_4$&$e_2 e_1=e_3+\af e_4$  \\
 $  \lb\ne 0, \af\ne 0$&& $e_2 e_2=e_3$ & $e_2e_3=\bt e_4$&$e_3e_3=e_4$\\
 
 \hline$\cd {4}{106}(\af)$&$:$& $e_1 e_3 = e_4$ & $e_2 e_1=e_3+\af e_4$  & $e_2 e_2=e_3$ &$e_3e_3=e_4$\\
 
 \hline$\cd {4}{107}(\lambda)$&$:$& 
 $e_1 e_1 = \lambda e_3$ & $e_1e_3=\Theta e_4$& $e_2 e_1=e_3$  \\
 && $e_2 e_2=e_3$ & $e_2e_3=e_4$&$e_3e_3=e_4$\\
 
 \hline$\cd {4}{108}(\lambda)$&$:$& 
 $e_1 e_1 = \lambda e_3$ & $e_1e_3=(1-\Theta) e_4$& $e_2 e_1=e_3$  \\
 $\lb \not\in \{0, \frac{1}{4}\}$&& $e_2 e_2=e_3$ & $e_2e_3=e_4$&$e_3e_3=e_4$\\
 
 \hline$\cd {4}{109}(\lambda,\af)$&$:$& $e_1 e_1 = \lambda e_3$ & $e_2 e_1=e_3+e_4$  \\&& $e_2 e_2=e_3$ &$e_2e_3=\alpha e_4$&$e_3e_3=e_4$\\

  \hline$\cd {4}{110}(\lambda)$&$:$& $e_1 e_1 = \lambda e_3$ & $e_2 e_1=e_3$  & $e_2 e_2=e_3$ \\&&$e_2e_3= e_4$&$e_3e_3=e_4$\\
 
   \hline$\cd {4}{111}(\lambda)$&$:$& $e_1 e_1 = \lambda e_3$ & $e_2 e_1=e_3$  & $e_2 e_2=e_3$ &$e_3e_3=e_4$\\
 
\hline$\cd {4}{112}(\lambda, \af, \bt, \gamma)$&$:$& 
$e_1 e_1 = \lambda e_3+e_4$ & $e_1e_3=\af e_4$ & $e_2 e_1=e_3+\bt e_4$  \\
&&$e_2 e_2=e_3$& $e_2e_3=\gamma e_4$&$e_3e_3=e_4$\\

\hline
$\n^4_{01}(\af)$ &$:$& 
$e_1 e_1 = e_2$& $e_1e_2=\af e_4$& $e_1e_3=e_4$ \\& & $e_2e_1=e_4$& $e_2 e_2=e_3$ \\ 
\hline
$\n^4_{02}(\af,\bt,\gm,\delta)$ &$:$& 
$e_1 e_1 = e_2$& $e_1e_2=\af e_4$& $e_1e_3=\gm e_4$&  $e_2e_1=\bt e_4$\\&& 
$e_2 e_2=e_3$& $e_2e_3=\delta e_4$& $e_3e_2=e_4$& $e_3e_3=e_4$\\ 
\hline
$\n^4_{03}(\af)$ &$:$& $e_1 e_1 = e_2$& $e_1e_2=\af e_4$& $e_1e_3=e_4$\\&& $e_2 e_2=e_3$ & $e_2e_3=e_4$\\  \hline

$\n^4_{04}(\af,\bt,\gm)$ &$:$& 
$e_1 e_1 = e_2$&$e_1e_2=\af e_4$&$e_1e_3=\gm e_4$&$e_2e_1=\bt e_4$\\
&&  $e_2 e_2=e_3$&$e_2e_3=e_4$&$e_3e_3=e_4$ \\  \hline

$\n^4_{05}(\af,\bt)$ &$:$& 
$e_1 e_1 = e_2$&$e_1e_2=\af e_4$&$e_1e_3=e_4$ \\&&$e_2e_1=\bt e_4$&  $e_2 e_2=e_3$&$e_3e_3=e_4$ \\  \hline

$\n^4_{06}(\af)$ &$:$& $e_1 e_1 = e_2$&$e_1e_2=\af e_4$&$e_2e_1=e_4$\\&&  $e_2 e_2=e_3$&$e_3e_3=e_4$ \\  \hline

$\n^4_{07}$ &$:$& $e_1 e_1 = e_2$&$e_1e_2=e_4$&$e_1e_3=e_4$&  $e_2 e_2=e_3$ \\  \hline

$\n^4_{08}(\af,\bt)$ &$:$& $e_1 e_1 = e_2$&$e_1e_2=e_4$&$e_1e_3=\bt e_4$\\&&$e_2e_1=\af e_4$&  $e_2 e_2=e_3$ &$e_3e_1=e_4$\\  \hline

$\n^4_{09}(\af,\bt)$ &$:$& 
$e_1 e_1 = e_2$&$e_1e_2=\af e_4$&$e_1e_3=\bt e_4$ \\&&  $e_2 e_2=e_3$ &$e_2e_3=e_4$&$e_3e_1=e_4$ \\  \hline

$\n^4_{10}$ &$:$& $e_1 e_1 = e_2$&$e_1e_2=e_4$&  $e_2 e_2=e_3$ &$e_2e_3=e_4$ \\  \hline

$\n^4_{11}$ &$:$& $e_1 e_1 = e_2$&$e_1e_2=e_4$&  $e_2 e_2=e_3$&$e_3e_3=e_4$ \\  \hline

$\n^4_{12}(\af)$ &$:$& $e_1 e_1 = e_2$&  $e_1e_3=\af e_4$&$e_2e_1=e_4$\\& &  $e_2 e_2=e_3$ &$e_3e_1=e_4$ \\  \hline

$\n^4_{13}(\af,\bt,\gm)$ &$:$& 
$e_1 e_1 = e_2$&$e_1e_3=\bt e_4$&$e_2e_1=\af e_4$&  $e_2 e_2=e_3$\\
&&$e_2e_3=\gm e_4$&$e_3e_1=e_4$&$e_3e_2=e_4$ \\  \hline

$\n^4_{14}(\af, \bt)$ &$:$& $e_1 e_1 = e_2$&$e_1e_3 =e_4$&$e_2e_1= \af e_4$ \\&&  $e_2 e_2=e_3$ &$e_2e_3=\bt e_4$&$e_3e_2=e_4$ \\  \hline

$\n^4_{15}(\af)$ &$:$& $e_1 e_1 = e_2$&$e_2e_1=e_4$&  $e_2 e_2=e_3$ \\&&$e_2e_3=\af e_4$ &$e_3e_2=e_4$ \\  \hline

$\n^4_{16}(\af)$ &$:$& $e_1 e_1 = e_2$&$e_1e_3=\af e_4$ &  $e_2 e_2=e_3$ &$e_3e_1=e_4$\\  \hline

$\n^4_{17}$ &$:$& $e_1 e_1 = e_2$&  $e_1e_3=e_4$ & $e_2 e_2=e_3$\\  \hline

$\n^4_{18}$ &$:$& $e_1 e_1 = e_2$&  $e_2 e_2=e_3$ &$e_2e_3=e_4$\\  \hline

$\n^4_{19}(\af)$ &$:$& $e_1 e_1 = e_2$&  $e_2 e_2=e_3$ &$e_2e_3=\af e_4$ &$e_3e_2=e_4$ \\  \hline

$\n^4_{20}$ &$:$& $e_1 e_1 = e_2$&  $e_2 e_2=e_3$ &$e_3e_3=e_4$\\  \hline

$\n^4_{21}(\af,\bt)$ &$:$& $e_1 e_1 = e_2$&$e_1e_2=\af e_4$&$e_1e_3=e_4$ \\&&$e_2e_1=e_3+\bt e_4$&  $e_2 e_2=e_3$ \\  \hline

$\n^4_{22}(\af,\bt,\gm)$ &$:$& $e_1 e_1 = e_2$&$e_1e_2=\af e_4$&$e_1e_3=\gm e_4$ \\&&$e_2e_1=e_3+\bt e_4$&  $e_2 e_2=e_3$&$e_3e_1=e_4$ \\  \hline

$\n^4_{23}(\af,\bt,\gm)$ &$:$& 
$e_1 e_1 = e_2$&$e_1e_2=\af e_4$&$e_1e_3=\bt e_4$&$e_2e_1=e_3$\\
&&  $e_2 e_2=e_3$&$e_2e_3=e_4$&$e_3e_1=\gm e_4$ \\  \hline

$\n^4_{24}(\af,\bt,\gm,\delta)$ &$:$& 
$e_1 e_1 = e_2$&$e_1e_3=\bt e_4$&$e_2e_1=e_3+\af e_4$&  $e_2 e_2=e_3$ \\
&&$e_2e_3=\delta e_4$&$e_3e_1=\gm e_4$ & $e_3e_2=e_4$\\  \hline

$\n^4_{25}(\af, \bt, \gm, \delta, \eps)$ &$:$& 
$e_1 e_1 = e_2$&$e_1e_2= \af e_4$& $e_2e_1=e_3+\bt e_4$ & $e_2 e_2=e_3$\\
&&$e_2e_3=\delta e_4$ & $e_3e_1=\gm e_4$ & $e_3e_2=\eps e_4$ &$e_3e_3=e_4$\\  \hline

$\n^4_{26}(\af)$ &$:$& 
$e_1 e_1 = e_2$&$e_1e_2=e_4$&$e_1e_3=-\af e_4$ \\&&$e_2e_1=e_3$  &$e_2e_2=\af e_4$&$ e_3e_1=(1+2 \af) e_4$  \\  \hline
 
$\n^4_{27}(\af,\bt)$ &$:$& $e_1 e_1 = e_2$&$e_1e_2=\af e_4$&$e_2e_1=e_3$\\& &$e_2e_2=\bt e_4$&$e_2e_3=e_4$&$e_3e_1=e_4$&  \\  \hline

$\n^4_{28}(\af,\bt,\gm)$ &$:$& $e_1 e_1 = e_2$&$e_1e_2=\af e_4$&$e_2e_1=e_3$&$e_2e_2=\bt e_4$\\&&$e_2e_3=e_4$&$e_3e_1=\gm e_4$&$e_3e_3=e_4$&  \\  \hline

$\n^4_{29}(\af,\bt)$ &$:$& $e_1 e_1 = e_2$&$e_1e_2=\af e_4$&$e_2e_2=\bt e_4$\\& &$e_2e_1=e_3$&$e_3e_1=e_4$&$e_3e_3=e_4$&  \\  \hline

$\n^4_{30}(\af) $ &$:$& $e_1 e_1 = e_2$&$e_1e_2=\af e_4$&$e_2e_1=e_3$\\&&$e_2e_2=e_4$&$e_2e_3=e_4$  \\  \hline

$\n^4_{31}(\af) $ &$:$& $e_1 e_1 = e_2$&$e_1e_2=\af e_4$&$e_2e_1=e_3$\\&&$e_2e_2=e_4$&$e_3e_3=e_4$  \\  \hline

$\n^4_{32}$ &$:$& $e_1 e_1 = e_2$&$e_1e_2=e_4$&$e_2e_1=e_3$&$e_2e_3=e_4$  \\  \hline

$\n^4_{33}$ &$:$& $e_1 e_1 = e_2$&$e_1e_2=e_4$&$e_2e_1=e_3$&$e_3e_3=e_4$  \\  \hline

$\n^4_{34}(\af,\bt)$ &$:$& 
$e_1 e_1 = e_2$&$e_1e_3=\bt e_4$&$e_2e_1=e_3$\\&&$e_2e_2=\af e_4$&
\multicolumn{2}{l}{$e_3e_1=(1-2\bt) e_4$}  \\  \hline

$\n^4_{35}(\af)$ &$:$& $e_1 e_1 = e_2$&$e_1e_3=\af e_4$&$e_2e_1=e_3$\\& &$e_2e_2=e_4$&$e_3e_1=-2\af e_4$&$e_3e_2=e_4$&  \\  \hline

$\n^4_{36}(\af,\bt)$ &$:$& 
$e_1 e_1 = e_2$&$e_2e_1=e_3$&$e_2e_2=\af e_4$\\
$\bt\neq 0$&  &$e_2e_3=\bt e_4$&$e_3e_1=e_4$&$e_3e_2=e_4$&  \\  \hline

$\n^4_{37}(\af)$ &$:$& 
$e_1 e_1 = e_2$&$e_2e_1=e_3$&$e_2e_2=e_4$ \\
$\af\neq0$&&$e_2e_3=\af e_4$&$e_3e_2=e_4$  \\  \hline

$\n^4_{38}$ &$:$& $e_1 e_1 = e_2$&$e_1e_3=e_4$&$e_2e_1=e_3$\\&&$e_3e_1=-2e_4$&$e_3e_2=e_4$ \\  \hline

$\n^4_{39}$ &$:$& $e_1 e_1 = e_2$&$e_2e_1=e_3$&$e_2e_3=e_4$&  \\  \hline

$\n^4_{40}(\af)$ &$:$& $e_1 e_1 = e_2$&$e_2e_1=e_3$&$e_2e_3=\af e_4$&$e_3e_2=e_4$&  \\  \hline

$\n^4_{41}$ &$:$& $e_1 e_1 = e_2$&$e_2e_1=e_3$&$e_3e_3=e_4$&  \\  \hline

$\n^4_{42}(\lb,\af)$ & $:$ & 
$e_1 e_1 = e_2$ & $e_1 e_2=e_3$ & \multicolumn{2}{l}{$e_1e_3=(\af(\lb-2)+1)e_4$}\\
$\lb\ne 1$&& $e_2 e_1=\lb e_3+e_4$ & 
\multicolumn{2}{l}{$e_2e_2=\frac{(\lb+1) ((\lb^2-1)\af + \lb + 2)}{1-\lb}e_4$} & $e_3e_1=(\af(1-2\lb)-2)e_4$  \\ 

\hline

$\n^4_{43}(\af,\bt)$ & $:$ & $e_1 e_1 = e_2$ & $e_1 e_2=e_3$ & $e_1e_3=-\frac {5\af}2e_4$ & $e_2 e_1=-\frac 12e_3+\bt e_4$\\
&& $e_2e_2=e_4$ & $e_2 e_3=-\frac 12 e_4$ & $e_3e_1=2\af e_4$ & $e_3e_2 = e_4$ \\

\hline

$\n^4_{44}(\af)$ & $:$ & $e_1 e_1 = e_2$ & $e_1 e_2=e_3$ & $e_1e_3=-\frac {5\af}2e_4$ & $e_2 e_1=-\frac 12e_3+e_4$\\
&& $e_2 e_3=-\frac 12 e_4$ & $e_3e_1=2\af e_4$ & $e_3e_2 = e_4$ \\

\hline

$\n^4_{45}(\lb,\af,\bt)$ & $:$ & $e_1 e_1 = e_2$ & $e_1 e_2=e_3$ & \multicolumn{2}{l}{$e_1e_3=(\af(\lb-2)+1)e_4$}\\
&& $e_2 e_1=\lb e_3$& $e_2e_2=\bt e_4$ & \multicolumn{2}{l}{$e_3e_1=(\af(1-2\lb)-2)e_4$}  \\ 

\hline

$\n^4_{46}(\lb,\af)$ & $:$ & $e_1 e_1 = e_2$ & $e_1 e_2=e_3$ & $e_1e_3=\af(\lb-2)e_4$ & $e_2 e_1=\lb e_3$\\
$\lb\ne-\frac 12$&& $e_2e_2=e_4$ &  $e_2 e_3=-\frac 12 e_4$ & $e_3e_1=\af(1-2\lb)e_4$ & $e_3e_2 = e_4$  \\ 

\hline

$\n^4_{47}$ & $:$ & $e_1 e_1 = e_2$ & $e_1 e_2=e_3$ & $e_1e_3=-\frac 52e_4$ & $e_2 e_1=-\frac 12e_3$\\
&& $e_2 e_3=-\frac 12 e_4$ & $e_3e_1=2 e_4$ & $e_3e_2 = e_4$ \\

\hline

$\n^4_{48}(\lb,\af,\bt)$ & $:$ & $e_1 e_1 = e_2$ & $e_1 e_2=e_3$ & $e_1e_3=e_4$ & $e_2 e_1=\lb e_3+\af e_4$\\
$\lb\ne -\frac 12$&& $e_2e_2=\bt e_4$ & $e_2 e_3=\lb e_4$ & $e_3e_1=-2e_4$ & $e_3e_2 = e_4$ \\

\hline

$\n^4_{49}(\af,\bt)$ & $:$ & $e_1 e_1 = e_2$ & $e_1 e_2=e_3$ & $e_1e_3=\bt e_4$\\
&& $e_2 e_1=-e_3+\af e_4$ & $e_2e_2=e_4$ & $e_2 e_3=e_4$\\
&& $e_3e_1=-2\bt e_4$ & $e_3e_2 = e_4$ & $e_3e_3=e_4$ \\

\hline

$\n^4_{50}(\lb,\af,\bt)$ & $:$ & $e_1 e_1 = e_2$ & $e_1 e_2=e_3$ & $e_1e_3=e_4$ & $e_2 e_1=\lb e_3+\af e_4$\\
$\lb\ne-1$&& $e_2e_2=\bt e_4$ & $e_3e_1=-2e_4$  & $e_3e_3=e_4$ \\

\hline

$\n^4_{51}(\lb,\af)$ & $:$ & $e_1 e_1 = e_2$ & $e_1 e_2=e_3$  & $e_2 e_1=\lb e_3+\af e_4$\\
$\lb\ne -\frac 12$&& $e_2e_2=e_4$  & $e_2e_3=\lb e_4$  & $e_3e_2=e_4$ \\

\hline

$\n^4_{52}(\lb,\af,\bt,\gm)$ & $:$ & $e_1 e_1 = e_2$ & $e_1 e_2=e_3$  & $e_2 e_1=\lb e_3+\af e_4$ & $e_2e_2=\bt e_4$\\
$\gm\ne 1$&& $e_2e_3=\gm e_4$  & $e_3e_2=e_4$ & $e_3e_3=e_4$ \\

\hline

$\n^4_{53}(\lb,\af,\bt)$ & $:$ & $e_1 e_1 = e_2$ & $e_1 e_2=e_3$  & $e_2 e_1=\lb e_3+\af e_4$\\
&& $e_2e_2=\bt e_4$ & $e_2e_3=e_4$  & $e_3e_3=e_4$ \\

\hline

$\n^4_{54}(\lb,\af)$ & $:$ & $e_1 e_1 = e_2$ & $e_1 e_2=e_3$  & $e_2 e_1=\lb e_3+\af e_4$\\
$\lb\ne -1$&& $e_2e_2=e_4$  & $e_3e_3=e_4$ \\

\hline

$\n^4_{55}(\lb)$ & $:$ & $e_1 e_1 = e_2$ & $e_1 e_2=e_3$  & $e_2 e_1=\lb e_3+e_4$\\
$\lb\ne -\frac 12$ && $e_2e_3=\lb e_4$  & $e_3e_2=e_4$ \\

\hline

$\n^4_{56}(\lb)$ & $:$ & $e_1 e_1 = e_2$ & $e_1 e_2=e_3$  & $e_2 e_1=\lb e_3+e_4$ & $e_3e_3=e_4$ \\

\hline 

$\n^4_{57}(\lb,\af)$ & $:$ & $e_1 e_1 = e_2$ & $e_1 e_2=e_3$ & $e_1e_3=e_4$ & $e_2 e_1=\lb e_3$\\
&& $e_2e_2=\af e_4$ & $e_2e_3=e_4$ & $e_3e_1=-2e_4$ \\

\hline

$\n^4_{58}(\lb,\af,\bt)$ & $:$ & $e_1 e_1 = e_2$ & $e_1 e_2=e_3$ & $e_1e_3=e_4$ & $e_2 e_1=\lb e_3$\\
$\bt\not\in\{-\frac 12,\lb\}$&& $e_2e_2=\af e_4$ & $e_2e_3=\bt e_4$ & $e_3e_1=-2e_4$ & $e_3e_2=e_4$ \\

\hline

$\n^4_{59}(\af,\bt)$ & $:$ & $e_1 e_1 = e_2$ & $e_1 e_2=e_3$ & $e_1e_3=\bt e_4$\\
$\af\ne 1$&& $e_2 e_1=-e_3$ & $e_2e_2=\af e_4$ & $e_2e_3=e_4$\\
&& $e_3e_1=-2\bt e_4$ & $e_3e_2=e_4$ & $e_3e_3=e_4$ \\

\hline

$\n^4_{60}(\af)$ & $:$ & $e_1 e_1 = e_2$ & $e_1 e_2=e_3$ & $e_1e_3=e_4$ & $e_2 e_1=-e_3$\\
$\af\ne 0$&& $e_2e_2=\af e_4$ & $e_3e_1=-2e_4$ & $e_3e_3=e_4$ \\

\hline

$\n^4_{61}(\lb)$ & $:$ & $e_1 e_1 = e_2$ & $e_1 e_2=e_3$ & $e_2 e_1=\lb e_3$\\
&& $e_2e_2=e_4$ & $e_2e_3=e_4$ \\

\hline

$\n^4_{62}(\lb,\af)$ & $:$ & $e_1 e_1 = e_2$ & $e_1 e_2=e_3$ & $e_2 e_1=\lb e_3$\\
$\af\not\in\{-\frac 12,\lb\}$&& $e_2e_2=e_4$ & $e_2e_3=\af e_4$ & $e_3e_2=e_4$ \\

\hline

$\n^4_{63}$ & $:$ & $e_1 e_1 = e_2$ & $e_1 e_2=e_3$ & $e_2 e_1=-e_3$\\
&& $e_2e_2=e_4$ & $e_3e_3=e_4$ \\

\hline

$\n^4_{64}(\lb)$ & $:$ & $e_1 e_1 = e_2$ & $e_1 e_2=e_3$ & $e_2 e_1=\lb e_3$ & $e_2e_3=e_4$ \\

\hline

$\n^4_{65}(\lb,\af)$ & $:$ & $e_1 e_1 = e_2$ & $e_1 e_2=e_3$ & $e_2 e_1=\lb e_3$\\
&& $e_2e_3=\af e_4$ & $e_3e_2=e_4$ \\

\hline

$\n^4_{66}(\lb)$ & $:$ & $e_1 e_1 = e_2$ & $e_1 e_2=e_3$ & $e_2 e_1=\lb e_3$ & $e_3e_3=e_4$\\

\end{longtable}}

All of these algebras are pairwise non-isomorphic, except for the following:

{\tiny \begin{longtable}{c}
$\D{4}{01}(\lambda,0,\beta) \cong \D{4}{02}(\lambda,0,\beta) \cong \D{4}{04}(\lambda,\beta),\quad \D{4}{01}(\lambda,\alpha,0)_{\alpha \neq -1} \cong \D{4}{02}(\lambda,\alpha,0) \cong \D{4}{10}(\lambda,\alpha),\quad \D{4}{01}(\lambda,-1,0) \cong \D{4}{11}(\lambda,0),$\\

$\D{4}{03}(\lambda,0) \cong \D{4}{09}(\lambda,0),\quad \D{4}{03}\left(\lambda,(1-\Theta)^{-1}\right)_{\lambda \neq 0} \cong \D{4}{05}(\lambda,0)_{\lambda \neq 0}, \D{4}{03}\left(\lambda,\Theta^{-1}\right)\cong \D{4}{06}(\lambda,0),\quad \D{4}{04}(\lambda,0) \cong \D{4}{10}(\lambda,0),$\\

$\D{4}{05}(1/4,\alpha) \cong \D{4}{06}(1/4,\alpha),\quad \D{4}{07}(1/4) \cong \D{4}{08}(1/4), \quad
\D{4}{05}(0,\alpha) \cong \D{4}{07}(0) \cong \D{4}{23}(0) \cong  \D{4}{25}(0) \cong  \D{4}{40}(0),$\\

$\D{4}{12}(\lambda,0) \cong \D{4}{18}(\lambda,0),\quad  \D{4}{12}(1/4,\alpha) \cong \D{4}{13}(1/4,\alpha),\quad \D{4}{12}(0,\alpha)_{\alpha \neq -1} \cong \D{4}{14}(0,\alpha),\quad \D{4}{12}(0,-1) \cong \D{4}{17}(0),$\\

$\D{4}{13}(\lambda,0) \cong \D{4}{19}(\lambda,0),\quad \D{4}{14}(\lambda,0) \cong \D{4}{20}(\lambda,0),\quad \D{4}{14}(1/4,\alpha) \cong \D{4}{15}(1/4,\alpha),\quad \D{4}{15}(\lambda,0) \cong \D{4}{21}(\lambda,0),$\\

$\D{4}{18}(1/4,\alpha) \cong \D{4}{19}(1/4,\alpha),\quad \D{4}{18}(0,0) \cong \D{4}{22}(0) \cong  \D{4}{24}(0),\quad \D{4}{18}(1/4,-1) \cong \D{4}{19}(1/4,-1) \cong \D{4}{30}(1/4) \cong \D{4}{31}(1/4),$\\

$\D{4}{20}(1/4,\alpha) \cong \D{4}{21}(1/4,\alpha),\quad \D{4}{20}(1/4,-1) \cong \D{4}{21}(1/4,-1) \cong \D{4}{32}(1/4) \cong \D{4}{33}(1/4),$\\

$\D{4}{22}(1/4) \cong \D{4}{23}(1/4) \cong \D{4}{24}(1/4) \cong \D{4}{25}(1/4) \cong \D{4}{26}(1/4) \cong \D{4}{27}(1/4) \cong \D{4}{28}(1/4) \cong \D{4}{29}(1/4),$\\

$ \D{4}{37}(1/4) \cong \D{4}{38}(1/4),\quad \D{4}{39}(1/4) \cong \D{4}{40}(1/4).$
\end{longtable}
\begin{longtable}{lll}
$\cd 4{43}(\af)\cong\cd 4{43}(-\af)$ &
$\cd 4{44}(\af,\bt,\gm)\cong\cd 4{44}(\af,-\bt,-\gm)$& 
$\cd 4{47}(\af,\bt)\cong \cd 4{47}(\af,-\bt)$\\
$\cd 4{50}(\af)=\cd 4{50}(-\af)$& 
$\cd 4{54}(\af)\cong\cd 4{54}(-\af-1)$&
$\cd 4{57}(\af,\bt)\cong \cd 4{57}(\af+\bt,-\bt)$\\
$\cd 4{59}(\af,\bt)\cong\cd 4{59}(\af,-\bt)$&
$\cd {4}{91}(\lb, \af)\cong\cd {4}{91}(\lb, -\af)$ &
$\cd {4}{92}(\lb, \af)\cong\cd {4}{92}(\lb, -\af)$ \\ 
$\cd {4}{93}(\af)\cong\cd {4}{93}(-\af)$ &
$\cd {4}{94}(\af,\bt)\cong\cd {4}{94}(-\af,\bt)$ &
$\cd {4}{95}(\af)\cong\cd {4}{95}(-\af)$ \\
$\cd {4}{100}(\af)\cong\cd {4}{100}(-\af)$ &
$\cd {4}{101}(\af,\bt)\cong\cd {4}{101}(-\af,-\bt)$ &
$\cd {4}{109}(\lb,\af)\cong\cd {4}{109}(\lb,-\af)$ \\
\multicolumn{3}{c}{$\cd {4}{112}(\lb,\af,\bt,\gm)\cong\cd {4}{112}(\lb,-\af,\bt,-\gm)$} \\
\multicolumn{3}{c}{$\cd {4}{112}(\lb,\af,\bt,\gm)\cong \cd {4}{112}\left(\lb,(\gm-\af\bt)\sqrt{\frac{-\lb}{1-\bt+\lb\bt^2}},\frac 1\lb-\bt,(\frac\gm\lb-\frac\af\lb-\bt\gm)\sqrt{\frac{-\lb}{1-\bt+\lb\bt^2}}\right)$, if $\lb\ne 0$, $\bt\ne\frac{1\pm\sqrt{1-4\lb}}{2\lb}$}\\

$\n^4_{02}(\af,\bt,\gm,\delta) \cong \n^4_{02}(-(\af,\bt,\gm),\delta)$ &
$\n^4_{04}(\af,\bt,\gm) \cong \n^4_{04}(-(\af,\bt,\gm))$ & $\n^4_{05}(\af,\bt) \cong \n^4_{05}(\sqrt[3]{1}(\af,\bt))$ \\
$\n^4_{29}(\af,\bt) \cong \n^4_{29}(-\af,\bt)$ &
$\n^4_{31}(\af) \cong \n^4_{31}(-\af)$ & $\n^4_{50}(\lb,\af,\bt) \cong \n^4_{50}(\lb,-\af,\bt)$ \\ &
$\n^4_{54}(\lb,\af) \cong \n^4_{54}(\lb,-\af)$

\end{longtable} 
} 
 
\end{theorem} 

\section{Applications}\label{S:apps}

\subsection{The algebraic classification of $4$-dimensional nilpotent Lie-admissible  algebras}
The variety of Lie-admissible algebras is defined by the following identity 
$$[[x,y],z]+[[y,z],x]+[[z,x],y]=0,$$
where $[x,y]=xy-yx.$
Lie-admissible algebras satisfy the following fundamental property: {\it under the commutator multiplication each Lie-admissible algebra is a Lie algebra}.

\begin{corollary}
Let $\bf A$ be a complex $4$-dimensional nilpotent Lie-admissible algebra.
Then $\bf A$ is isomorphic to an algebra from the list given in Theorem \ref{teo-alg}.
 \end{corollary}
\begin{Proof}
Thanks to \cite{kkl19} each $4$-dimensional nilpotent anticommutative algebra is Lie.
Hence, each $4$-dimensional nilpotent algebra is a Lie-admissible algebra.
\end{Proof}

\subsection{The algebraic classification of 
$4$-dimensional nilpotent Alia type algebras}
Let $\bf A$ be an algebra with a certain bilinear multiplication $(x,y)\mapsto xy$, which we call standard.
Let us define 
\begin{center}$\textsf{T}(x,y,z, {\textsf P})={\textsf P}([x,y],z)+{\textsf P}([y,z],x)+{\textsf P}([z,x],y)$
\end{center}
where $[x,y]=xy-yx$ and ${\textsf P}$ is a bilinear multiplication on the same vector space. 
Now we are ready to introduce Alia type algebras, which appeared in \cite{dzhuma1, dzhuma2}: 
\begin{enumerate}
  
   \item the variety of $0$-Alia (also known as $0$-anti-Lie-admissible) algebras is defined by the identity $\textsf{T}(x,y,z, {\textsf P})=0,$  where ${\textsf P}$ is the standard multiplication; 
     
   \item the variety of $1$-Alia (also known as $1$-anti-Lie-admissible) algebras is defined by the identity $\textsf{T}(x,y,z, {\textsf P})=0,$  where ${\textsf P}$ is the  Jordan commutator multiplication ${\textsf P}(x,y)=xy+yx;$ 
    
   \item the variety of two-sided Alia (also known as anti-Lie-admissible) algebras is defined by the identities $\textsf{T}(x,y,z, {\textsf P}_1)=0$ and $\textsf{T}(x,y,z, {\textsf P}_2)=0,$  where ${\textsf P}_1$ is the standard multiplication and  ${\textsf P}_2$ is the opposite  multiplication.

    \end{enumerate}

\begin{corollary}
Let $\bf A$ be a complex $4$-dimensional nilpotent  $0$-Alia ($1$-Alia, or two sided Alia)  algebra.
Then $\bf A$ is isomorphic to an algebra from the list given in Theorem \ref{teo-alg}, except

\begin{longtable}{c}
$\cd {4}{79} - \cd {4}{85}, \cd {4}{87} - \cd {4}{112}, \n^4_{25}, \  \n^4_{28}, \ \n^4_{29}, \ \n^4_{31}, \ \n^4_{33}, \ \n^4_{41}, \
\n^4_{49}, \ \n^4_{50}(\lambda\neq 1), $\\
$\ \n^4_{52}(\lambda\neq 1),\ \n^4_{53}(\lambda\neq 1), \ \n^4_{54}(\lambda\neq 1),\ \n^4_{56}(\lambda\neq 1),\n^4_{59}, \n^4_{60}, \n^4_{63}, \n^4_{66}(\lambda\neq 1).$
\end{longtable}

 \end{corollary}
\begin{Proof} 
After verification, we have that all $3$-dimensional nilpotent algebras are
$0$-Alia,  $1$-Alia, and two-sided Alia.
Analyzing the cocycles on $3$-dimensional algebras we see that the subspaces of $0$-Alia cocycles,  $1$-Alia cocycles, and two-sided Alia cocycles coincide and are listed below.
\begin{longtable}{|l|l|}
\hline
Algebra  &   Cocycles   \\
\hline
${\mathfrak{CD}}_{01}^{3}, {\mathfrak{CD}}_{04}^{3}(1),
{\mathfrak{CD}}_{01}^{3*},{\mathfrak{CD}}_{02}^{3*}$ &  
 $\langle \Delta_{ij} \rangle_{1\leq i,j\leq 3} $ \\ 
 
\hline

${\mathfrak{CD}}_{02}^{3}, {\mathfrak{CD}}_{03}^{3}, {\mathfrak{CD}}_{04}^{3}(\lambda \neq 1),
{\mathfrak{CD}}_{03}^{3*},{\mathfrak{CD}}_{04}^{3*}$ &  
 $\langle \Delta_{ij} \rangle_{1\leq i,j\leq 3; (i,j) \neq (3,3) } $ \\ 
 
\hline
\end{longtable}	
It remains to choose the algebras from Theorem \ref{teo-alg} that are determined by Alia cocycles.

\end{Proof}


\begin{thebibliography}{}
\bibitem{ack}
Abdelwahab H.,  Calder\'on A.J., Kaygorodov I.,
    The algebraic and geometric classification of nilpotent binary Lie algebras, 
     International Journal of Algebra and Computation, 29 (2019), 6, 1113--1129.
  
    
 
     
\bibitem{cfk18}
Calderón Martín A., Fern\'andez Ouaridi A., Kaygorodov I.,
    On the classification of bilinear maps with radical of a fixed codimension,
    Linear and Multilinear Algebra, 2020, DOI:10.1080/03081087.2020.1849001




\bibitem{lisa}
Camacho L., Kaygorodov I.,  Lopatkin V., Salim M., 
    The variety of dual Mock-Lie algebras,  
    Communications in Mathematics,    28  (2020), 2,  161--178. 

  
  




\bibitem{degr3}
Cicalò S., De Graaf W.,   Schneider C.,
    Six-dimensional nilpotent Lie algebras,
    Linear Algebra and its Applications, 436 (2012), 1, 163--189.



\bibitem{usefi1}
Darijani I., Usefi H.,
    The classification of 5-dimensional $p$-nilpotent restricted Lie algebras over perfect fields. I,
    Journal of Algebra, 464 (2016), 97--140.


\bibitem{degr2}
De Graaf W., 
    Classification of 6-dimensional nilpotent Lie algebras over fields of characteristic not $2$, 
    Journal of Algebra, 309  (2007), 2, 640--653.

\bibitem{degr1}
De Graaf W., 
    Classification of nilpotent associative algebras of small dimension,
    International Journal of Algebra and Computation, 28 (2018),  1, 133--161.


\bibitem{demir}
Demir I., Misra K.,  Stitzinger E.,
    On classification of four-dimensional nilpotent Leibniz algebras,
    Communications in Algebra, 45 (2017), 3, 1012--1018.

  

\bibitem{dzhuma1}
Dzhumadildaev A., Bakirova A., 
    Simple two-sided anti-Lie-admissible algebras, 
    Journal of Mathematical Sciences (N.Y.), 161 (2009),   1, 31--36

 
 \bibitem{dzhuma2}
Dzhumadildaev A., Tulenbaev K.,  
    Exceptional 0-Alia algebras, 
    Journal of Mathematical Sciences (N.Y.), 161 (2009), 1, 37--40.
    
\bibitem{fkkv}
Fern\'andez Ouaridi A.,  Kaygorodov I.,  Khrypchenko M., Volkov Yu., 
    Degenerations of nilpotent algebras,
    arXiv:1905.05361
 

\bibitem{ha16}
Hegazi A., Abdelwahab H.,
    Classification of five-dimensional nilpotent Jordan algebras,
    Linear Algebra and its Applications, 494 (2016), 165--218.

 

\bibitem{hac16}
Hegazi A., Abdelwahab H., Calderón Martín A.,
    The classification of $n$-dimensional non-Lie Malcev algebras with $(n-4)$-dimensional annihilator, 
    Linear Algebra and its Applications, 505 (2016), 32--56.

 
    
 
 
 

 
 
 
     
\bibitem{kk20}
Kaygorodov I., Khrypchenko M., 
    The algebraic  classification of nilpotent $\mathfrak{CD}$-algebras,
    Linear and Multilinear algebra, 2020, DOI: 10.1080/03081087.2020.1856030
    
 


\bibitem{kkl19}
Kaygorodov I., Khrypchenko M., Lopes S.,  
    The algebraic and geometric classification of nilpotent anticommutative algebras, 
    Journal of Pure and Applied Algebra,  224  (2020), 8, 106337. 
 
 
 
\bibitem{kkp19geo}
Kaygorodov I.,  Khrypchenko M.,  Popov Yu.,
    The algebraic and geometric classification of nilpotent terminal algebras,
    Journal of Pure and Applied Algebra,   225  (2021), 6, 106625.



\bibitem{omirov}
Kaygorodov I., Lopes S., Páez-Guillán P.,
    Non-associative central extensions of null-filiform associative algebras,
    Journal of Algebra, 560 (2020), 1190--1210.
  
  
  
 

\bibitem{kpv19}
Kaygorodov I., P\'{a}ez-Guill\'{a}n  P., Voronin V.,  
    The algebraic and geometric classification of nilpotent bicommutative algebras,
      Algebras and Representation Theory,  23  (2020), 6, 2331-2347.


  
 
 


\bibitem{maz79}
Mazzola G.,
    The algebraic and geometric classification of associative algebras of dimension five, 
    Manuscripta Mathematica, 27 (1979), 1, 81--101. 
    
 

 

\bibitem{ss78}
Skjelbred T., Sund T.,
    Sur la classification des algebres de Lie nilpotentes,
    C. R. Acad. Sci. Paris Ser. A-B, 286 (1978), 5,  A241--A242.

 



 
 
 
 
 

















\end{thebibliography}
\end{document}